\documentclass[3p]{elsarticle}

\usepackage[latin1]{inputenc}
\usepackage[english]{babel}
\usepackage{graphicx}
\usepackage{epsfig}
\usepackage{color}
\usepackage{amsmath,amssymb}
\usepackage{amsthm}
\usepackage{bm}
\usepackage[normalem]{ulem}
\usepackage{caption}
\usepackage{subcaption}
\usepackage{url}
\usepackage{multirow}

\newcommand{\head}[1]{\textnormal{\textbf{#1}}}
\usepackage{nuovicomandi}
\newtheorem*{remark}{Remark}

\begin{document}


\begin{frontmatter}
\title{Uncertainty Quantification of geochemical and mechanical compaction in layered sedimentary basins}

\author[DICA]{Ivo Colombo}\ead{ivo.colombo@polimi.it}

\author[EPFL]{Fabio Nobile}\ead{ fabio.nobile@epfl.ch}

\author[DICA]{Giovanni Porta}\ead{giovanni.porta@polimi.it}

\author[MOX]{Anna Scotti}\ead{anna.scotti@polimi.it}

\author[IMATI]{Lorenzo Tamellini\corref{cor1}}\ead{tamellini@imati.cnr.it}

\cortext[cor1]{Corresponding author} 

\address[DICA]{Dipartimento di Ingegneria Civile e Ambientale, Politecnico di Milano, Piazza Leonardo da Vinci 32, 20133 Milano, Italy}
\address[EPFL]{CSQI - MATHICSE, \'Ecole Polytechnique F\'ed\'erale de Lausanne, Station 8, CH 1015, Lausanne, Switzerland}
\address[MOX]{MOX, Dipartimento di Matematica, Politecnico di Milano, Via Bonardi 9, 20133  Milano, Italy}
\address[IMATI]{Istituto di Matematica Applicata e Tecnologie Informatiche ``E. Magenes'' - CNR, Via Ferrata 5, 27100 Pavia, Italy}

\begin{abstract}
In this work we propose an Uncertainty Quantification methodology for sedimentary basins evolution
under mechanical and geochemical compaction processes, which we model as a coupled,
time-dependent, non-linear, monodimensional (depth-only) system of PDEs with uncertain parameters.
While in previous works (Formaggia et al. 2013, Porta et al., 2014) 
we assumed a simplified depositional history with only one material, in this work
we consider multi-layered basins, in which each layer is characterized
by a different material, and hence by different properties.
This setting requires several improvements with respect to our earlier works,
both concerning the deterministic solver and the stochastic discretization.
On the deterministic side, we replace the previous fixed-point iterative
solver with a more efficient Newton solver at each step of the time-discretization.
On the stochastic side, the multi-layered structure
gives rise to discontinuities in the dependence of the state variables on
the uncertain parameters, that need an appropriate treatment for
surrogate modeling techniques, such as sparse grids, to be effective. We propose an innovative methodology to this end which relies on a change of coordinate system to align the discontinuities of the target function within the random parameter space. The reference coordinate system is built upon exploiting physical features of the problem at hand. We employ the locations of material interfaces, which display a smooth dependence on the random parameters and are therefore amenable to sparse grid polynomial approximations. We showcase the capabilities of our numerical methodologies through two synthetic test cases. In particular, we show that our methodology reproduces with high accuracy multi-modal probability density functions displayed by target state variables (e.g., porosity).
\end{abstract}


\begin{keyword}
Sedimentary basin modeling \sep
Uncertainty Quantification \sep
Random PDEs \sep
Sparse grids \sep
Stochastic Collocation Method 

\MSC[2010]
  41A10 \sep 
  65C20 \sep 
  65N30 
\end{keyword}

\end{frontmatter}


\section{Introduction}

Sedimentary basins occupy depressions of the Earth's crust where different materials may deposit along geologic times. Numerical simulation of compaction processes in sedimentary basins is relevant to a number of fields, e.g. for the characterization of petroleum systems in terms of hydrocarbon generation and migration \cite{TuncayOrtoleva04, Taylor2010}, understanding of large scale hydrologic behavior (e.g., compaction-driven flow or development of fluid overpressure)  \cite{Mello1994, Hunt1998, JiaoZheng1998, Osborne1999, McPherson2001, Neuzil1994, Neuzil2003, MarinMoreno2013, Colombo16_HJ}, or modeling the formation of ore deposits \cite{Wieck1995}.

Modeling basin scale compaction requires to consider mechanical compaction due to the sediments overburden. Moreover, geochemical processes may also occur, due to chemical precipitation and/or dissolution of minerals, and these may heavily affect the effective properties of the system (e.g., porosity, permeability). The mechanics and fluid dynamics of the system evolve as a result of these coupled processes \cite{Bethke, Wangen2010}.
Key target outputs of basin compaction models are the porosity, pressure and temperature spatial distributions along the basin history. 

The geomechanical evolution of subsurface systems and the characterization of fluid dynamics in geological media are classical applications for Uncertainty Quantification methodologies. This can be explained upon considering that our knowledge of \Lc{the architecture}{geological and physical aspects} and of the multi-scale spatial heterogeneity of the subsurface is typically incomplete. In the case of basin compaction we also deal with large characteristic evolutionary scales (millions of years, Ma) and considerable spatial dimensions (km). In this context direct measurements for the characterization of the key processes at the pertinent scale are typically scarce if not lacking altogether, and therefore the boundary conditions and the model parameters are generally poorly constrained. A common approach to deal with this lack of knowledge is to model the uncertain parameters as random variables, and to consider the model predictions as the outputs of a random input-output map, to be analyzed with statistical techniques.

In previous contributions we have developed a forward and inverse modeling technique for basin scale compaction under uncertainty \cite{feal:compgeo,lever.eal:inversion}, in a simplified framework in which we assume that compaction mainly takes place along the vertical direction. This allows for a relevant simplification of the model structure which is reduced to a one-dimensional (vertical) setting. In this context, our uncertainty quantification technique relies on a sparse grids polynomial surrogate model of the input-output map 
\cite{babuska.nobile.eal:stochastic2,barthelmann.novak.ritter:high,b.griebel:acta,nobile.tempone.eal:aniso,xiu.hesthaven:high}, which approximates the full compaction model outputs (e.g., porosity, temperature or pressure vertical profiles). 
The sparse grid construction requires first to solve the full compaction model for a number of parameter combinations, which is performed in \cite{feal:compgeo,lever.eal:inversion} through a fixed iteration algorithm. 
The sparse grid approximation can then be algebraically post-processed to obtain relevant information for the Uncertainty Quantification analysis, such as statistical moments and Sobol indices of the quantities of interest. Results by \cite{feal:compgeo,lever.eal:inversion} demonstrate that this approach is very effective for the quantification of uncertainty when the basin has sediments with homogeneous properties, as in such a case the key outputs are typically smooth functions of the model parameters. 

In this contribution we extend the approach to the case of a multi-layered material. Crucially, in the multi-layer case some 
state variables experience discontinuities across the different layers, since they typically assume different values depending on the properties associated with the geo-material. For example, the evolution of porosity can considerably change between sandstone and shale units, due to (i) the different mechanical properties associated to the sediments and (ii) geochemical processes leading to porosity changes which selectively take place as
a function of the local sediment composition (our model e.g. considers porosity reduction due to quartz cementation which
may largely affect porosity only in sandstone layers). The occurrence of such sharp variations across interfaces poses two key issues to the numerical approach proposed in \cite{feal:compgeo,lever.eal:inversion}. 

First, when stratified materials are of concern, permeability contrast of several order of magnitude are commonly found, e.g., across interfaces between permeable sandy layers and shale units characterized by very low permeability. For such low values, the fixed iteration method proposed in \cite{feal:compgeo} is prone to failure. In this work we solve this issue by implementing a Newton iteration algorithm, which considers the fully coupled system of mass, momentum and energy conservation along the basin depth, together with the constitutive relations which are associated with geochemical reactions .

A second key challenge of multi-layered cases is that when multiple model parameters are considered to be random, a single space-time location may be associated with the presence of different geo-materials when different realizations of the random parameters are considered. This results in a discontinuous input-output map and a sparse grid (global) polynomial approximation of it will typically feature a loss of accuracy and fairly low convergence rate.

The occurrence of discontinuities or sharp gradients in input-output maps is rather frequent in mathematical models describing physical problems which are of interest in the context of engineering and applied mechanics, such as conservation laws or advection-dominated transport processes. Special numerical methodologies and procedures have been recently proposed to address this issue in various contexts. 
Discontinuities in the parameters space have been typically tackled in literature with multi-element approximations \cite{Barth16_CCE,tryoen.lem.eal:intrusive_gal,foo.karniadakis:me-pcm2,wan.karniadakis:adaptive}, multi-wavelet approaches \cite{Tryoen2012,LeMaitre1,LeMaitre2}, discontinuities-tracking algorithm \cite{jakeman.archi.xiu:discont,jakeman.eal:minimal-discont,Gorodetsky.marzouk:discontinuities,webster:discontinuities,Lang15_CM,LiStinis15_JCP,sargsyan.eal:discontinuities},
or by suitably choosing or enriching the polynomial approximation basis \cite{PetterssonTchelepi16_CMAME,chantrasmi.eal:discontinous}.
Another possibility is to transform the target variables into proxies which are more amenable to polynomial approximation, 
see e.g. \cite{alexanderian1,alexanderian2,liao.zhang:transform}.

We propose here an alternative approach based on the observation that, although the state variables at a fixed depth depend in a discontinuous fashion on the uncertain parameters, the depth of the interfaces among layers actually features a smooth dependence on the parameters and can therefore be accurately approximated by a sparse-grids approach. Upon estimating the interface depth location, we can then map each realization onto a reference domain in which the discontinuities in depth are aligned, and perform a second sparse-grids interpolation in this reference coordinate system, in which the layers never mix and hence sparse  grids approximation of the state variables is effective. 

While our approach is clearly different from multi-element approximations in the parameter domain, we remark that  
it is not a discontinuity-tracking algorithm in the parameter space either; in other words, we stress the fact that 
we do not try to approximate the boundary of the regions of the parameter space that generate
porosity profiles (or any other state variables of interest) such that a specific geomaterial will be found
at a specific depth. This approach would actually be rather unconvenient and quite challenging in high dimensional parameter spaces, given that such region 
will in general depend on the depth of interest, therefore we would need to track not one but multiple discontinuities
in the parameter space. Instead,
for each value of the parameter, we predict the location of the \emph{physical discontinuities}, that can be accurately predicted
by a sparse grid by exploiting the physical properties of our specific problem of interest.  
The rest of this work is organized as follows.
First, the mathematical model and its deterministic numerical
discretization are presented in Section \ref{sec:mod+discr}. 
We then describe the methodology employed to deal with discontinuities in the
parameters space in Section \ref{sec:met}. In particular,
we first recall some basics of sparse grid surrogate model construction 
and then motivate why discontinuities occur and propose a suitable numerical treatment. 
In Section \ref{application} we showcase the capabilites of our 
proposed approaches through synthetic test cases, characterized by a realistic evolution.  
Finally, some conclusions are presented in Section \ref{conclusions}.

\section{Basin Compaction Model}\label{sec:mod+discr}

In this section we provide a description of the mathematical formulation that we employ to model compaction phenomena at sedimentary basin scale and the corresponding numerical discretization strategy adopted.
We consider two main drivers for compaction: 
\begin{itemize}
\item \textbf{Mechanical compaction}, acting in the full rock domain and essentially due to the load of overburden sediments, is described as a rearrangement of the deposited grains and the ensuing reduction of porosity. Its effect is considered proportional to the effective stress which is defined as the difference between total stress (lithological pressure) and pore water pressure. 
\item \textbf{Geochemical compaction}, acting in the sand-rich materials and related to quartz dissolution and precipitation, provokes reduction of porosity once the mechanism starts: the interstices are increasingly filled with quartz crystals. The minimum values that porosity could reach in this way are significantly lower than those related to pure mechanical compaction. 
\end{itemize}

Following previous approaches \cite{Colombo16_HJ,Bethke,feal:compgeo}, we rely on the assumptions that (a) the most relevant phenomena take place mainly along the vertical direction and (b) the rock domain is assumed to be fully saturated by a single fluid characterized by uniform properties. Note that the first assumption enables us to consider a one-dimensional system described by the domain $\Omega$(t) = [$z_{bot}$(t), $z_{top}$(t)], where $z_{bot}$ and $z_{top}$ denote the bottom and the top of the rock domain respectively, which can both vary with time. These assumptions introduce considerable simplifications with respect to the reality of sedimentary systems, e.g. by overlooking the occurrence of density driven or multi-phase flows. At the same time our simplified approach is still capable of interpreting porosity and pressure distributions observed in real sedimentary systems \cite{Colombo16_HJ}. 
In this context a key advantage of introducing a simplified model is that uncertainty propagation can be thoroughly characterized at limited computational cost, a task which would be unaffordable if we were to explicitly model the evolution of a three-dimensional system along geological scales. 

\subsection{Mathematical model}\label{sec:model}

We describe the time evolution of sedimentary materials along depth through a coupled system of partial differential equations.
We impose mass conservation for the solid and fluid phase along the vertical direction by a standard formulation of the conservation laws
\begin{equation} \label{eq:solid-mass}
	\frac{\partial[(1-\phi) \rho^s]}{\partial t}+\frac{\partial[(1-\phi) \rho^s u^s]}{\partial z} = q^s
\end{equation}
\begin{equation} \label{eq:fluid-mass}
	\frac{\partial(\phi \rho^l)}{\partial t}+\frac{\partial(\phi \rho^l u^l )}{\partial z} = q^l
\end{equation}
where {$\phi$ [-] is the porosity of sediments; $u^i$ [m $\mbox{s}^{-1}$], $\rho^i$ [kg $\mbox{m}^{-3}$] and $q^i$ [kg $\mbox{m}^{-3}$ $\mbox{s}^{-1}$] are the velocity, the density and the source term of phase \textit{i} (where \textit{i} = \textit{s} or \textit{l},  for the solid phase or the fluid phase) respectively. The last term conveys the effect of solid and fluid generation, generally related to geochemical processes.

Displacement of pore water is defined as the difference between the velocity of fluid and solid phase and is described through the equation
\begin{equation} \label{eq:darcy}
	u^D = \phi (u^l-u^s) = \frac{K}{\mu^l} \left( \frac{\partial P}{\partial z} - \rho^l g \right)
\end{equation}
where $u^D$ [m $\mbox{s}^{-1}$] is the Darcy flux, $K$ [$\mbox{m}^2$] is the permeability, $\mu^l$ [kg $\mbox{m}^{-1}$ $\mbox{s}^{-1}$] is the fluid dynamic viscosity, $P$ [kg $\mbox{m}^{-1}$ $\mbox{s}^{-2}$] is the pore pressure, $g$ [m $\mbox{s}^{-2}$] is the gravitational acceleration. 

Following \cite{Wangen2010}
, we set
\begin{equation} \label{eq:k-phi}
	\log	K(\phi) = k_1 \phi - k_2 -15
\end{equation}
where the permeability $K$ is given in [$\mbox{m}^2$] and $k_1$ and $k_2$ are empirical fitting parameters. These two parameters are usually estimated considering experiments at laboratory scale (\cite{Wangen2010}) or analysis on core data (\cite{Nelson1994, Neuzil1994}). Even if this empirical formulation was developed for sandstone and sand-rich material, it is also widely employed for clays and shales \cite{Wangen2010, Nelson1994, Neuzil1994},  upon proper tuning of the values of $k_1$ and $k_2$. This modeling choice is consistent with 
published dataset  (see, e.g., \cite{Neuzil1994}). 

The effect of mechanical compaction is described following the approach of \cite{Schneider1994}
\begin{equation} \label{eq:mec-comp}
	\frac{d\phi_M}{dt}=-\beta (\phi_0 - \phi_f) \exp(-\beta \sigma_c) \frac{d \sigma_c}{dt}
\end{equation}
where
\begin{displaymath}
	\frac{d \cdot}{dt} = \frac{\partial \cdot}{\partial t} + u^s \frac{\partial \cdot}{\partial z}
\end{displaymath}
is the material derivative. $\frac{d\phi_M}{dt}$ is the material porosity variation due to mechanical compaction, $\beta$ [$\mbox{kg}^{-1}$ m $\mbox{s}^2$] is the porous medium (uniaxial) vertical compressibility, $\phi_0$ is the sediment porosity value at time of deposition, $\phi_f$ is the minimum porosity value achieved by pure mechanical compaction, and $\sigma_c$ [kg $\mbox{m}^{-1}$ $\mbox{s}^{-2}$] is the vertical effective stress, defined as the difference between lithostatic pressure and pore water pressure.

Quartz precipitation is considered modeled following \cite{Walderhaug1996}
\begin{equation} \label{eq:quartz}
	\frac{d \phi_Q}{dt} = A \frac{M_Q}{\rho_Q} a_q 10^{b_q T}  
\end{equation}
\begin{displaymath}
	A = A_0 \frac{\phi}{\phi_{act}}
\end{displaymath}
where $\phi_Q$ [-] is the quartz volumetric fraction; $M_Q$ [kg $\mbox{mol}^{-1}$] is the quartz molar mass; $\rho_Q$ [kg $\mbox{m}^{-3}$] is the density of quartz; $A_0$ [$\mbox{m}^{-1}$] and $\phi_{act}$ [-] are respectively the specific surface and porosity at the onset of quartz cementation; $a_q$ [mol $\mbox{m}^{-2}$ $\mbox{s}^{-1}$] and $b_q$ [$\mbox{K}^{-1}$] are system parameters. The activation of this reaction is temperature dependent: quartz starts to precipitate when temperature \textit{T} reaches a critical value $T_C$ that is generally assumed \Lc{to be variable between 353 and}{to range from 353 to} 373 K (\cite{Walderhaug1994}).

As a consequence of equations \eqref{eq:mec-comp}-\eqref{eq:quartz}, the porosity variation is then described by $d\phi = d\phi_M - d\phi_Q$, where $d\phi_M$ and $d\phi_Q$ denote the porosity variation due to mechanical compaction and quartz precipitation, respectively.

Temperature evolution inside the rock domain is modeled by 
\begin{equation} \label{eq:thermal}
	C_T \frac{dT}{dt} + \rho^l c^l u^D \frac{\partial T}{\partial z} - \frac{\partial}{\partial z} \left( K_T \frac{\partial T}{\partial z} \right)=0
\end{equation}
\begin{displaymath}
	C_T = \phi \rho^l c^l +  (1 - \phi) \rho^s c^s
\end{displaymath}
\begin{displaymath}	
	K_T = {\lambda_l}^{\phi} {\lambda_s}^{1-\phi}
\end{displaymath}

where $C_T(\phi)$ [kg $\mbox{m}^{-1}$ $\mbox{K}^{-1}$ $\mbox{s}^{-2}$] is the effective thermal capacity of the medium, $c^l$ and $c^s$ [$\mbox{m}^{2}$ $\mbox{K}^{-1}$ $\mbox{s}^{-2}$] are the liquid and solid specific thermal capacities, respectively; $K_T({\phi})$ [kg m $\mbox{K}^{-1}$ $\mbox{s}^{-3}$] is the thermal conductivity; and $\lambda_l$ and $\lambda_s$ are fluid and solid specific conductivities (\cite{Wangen2010}). Equation (\ref{eq:thermal}) includes dynamics of heat exchanges due to fluid advection, solid displacement and thermal diffusion. 

An appropriate selection of initial and boundary conditions is needed for the solution of the nonlinear system of partial differential equations (\ref{eq:fluid-mass})-(\ref{eq:thermal}). The basement that lies below the sedimentary basin is typically composed of rocks of igneous and metamorphic origin. We assume here that the basement can be considered as an impermeable layer, upon imposing zero fluid flux at $z  = z_{bot}$. We assume that the considered sedimentary system is submerged below the sea level, to which we assign constant depth $h_{sea}$, i.e. we set $z_{top} = -h_{sea}$. To take into account the presence of the sea water layer a fixed load proportional to the overlying water depth is assumed on the top of the basin, i.e. $p(z = z_{top}) = \gamma_{sea}h_{sea}$.  At the top location $z_{top}$ temperature is also assigned equal to $T(z = z_{top}) = T_{top}$ as a Dirichlet boundary condition. A heat flux is assigned as a Neumann boundary condition at the bottom of the basin, by fixing the local vertical derivative $\left. \frac{\partial T}{\partial z}\right|_{z = z_{bot}} = G_T$.





\subsection{Numerical discretization}\label{section:discretization}

The numerical discretization follows to a large extent the strategy described in \cite{feal:compgeo} concerning the space and time discretization of the equations, but differs significantly in the non-linear solver.

As in \cite{feal:compgeo} a Lagrangian approach  is
adopted to address the temporal evolution of the computational
domain, i.e. the computational grid is deformed under
the effect of compaction according to the solid matrix movement. We have, at each time step $t^k$, $N^k$ cells and $N^k+1$ nodes whose position is denoted as $z^{t^k}_i$ for $i=1,\ldots N^k+1$. Due to compaction the size of the elements of the grid, $h_i=z^{t^k}_{i+1}-z^{t^k}_i$, can change in time.
To include the deposition of new sediment layers over time we take into account sedimentation by
a modified load at the top of the basin until the thickness of fresh sediments equals the
characteristic size of the mesh elements. At that point, a new element
is added to the computational grid, see figure \ref{pic:griglia}.
\begin{figure}\centering
 \includegraphics[width=0.5\textwidth]{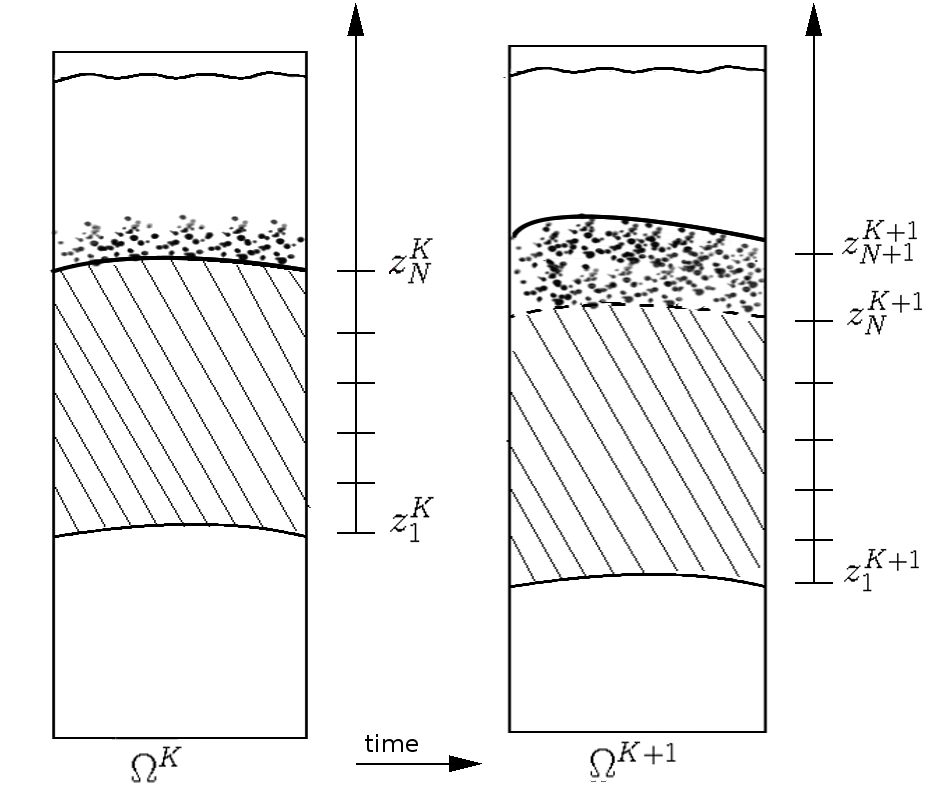}
 \caption{Deformation of the grid and addition of new elements with sedimentation}\label{pic:griglia}
\end{figure}

Thanks to the Lagrangian approach it is straightforward to account for the presence of different layers, characterized by different mechanical and chemical properties as required by the applications considered in this study. 
The evolution equations of the model are discretized in time with the Implicit Euler method. The size of the time step is chosen according to the averege sedimentation rate $U_\mathrm{sed}$ of the case of interest in such a way that a new cell grid of size $h$ is added at the top of the domain after $\alpha$ time step:

$$\delta t=\dfrac{h}{\alpha U_\mathrm{sed}}.$$

Concerning space discretization, mixed finite elements are used for the pressure and temperature problems, in particular pressure and temperature are approximated with piecewise constant $\mathbb{P}_0$ elements, while for Darcy velocity and heat flux we employ, since we are considering a 1D domain, simply nodal $\mathbb{P}_1$ elements.

\Lc{
For the numerical solution of the system of coupled equations previous works, such as \cite{feal:compgeo}, used a strategy inspired by the work of \cite{Mello2009}
related to three-dimensional basin modeling and based on an idea first presented in \cite{Chen}.
However that formulation of the iterative splitting has proven to be unstable for low permeable materials, due to the strong coupling between pressure and porosity. }
{
For the numerical solution of the resulting system of coupled equations a common approach, used in previous works such as  \cite{feal:compgeo}, consists in an iterative splitting, i.e. the problem is split into sub-problems that are solved in sequence, possibly iterating until convergence is achieved. However this strategy, in the formulation proposed by \cite{Mello2009} and based on an idea first presented in \cite{Chen},  has proven to be unstable for low permeable materials, due to the strong coupling between pressure and porosity.
}

For this reason, since we are interested in simulations that involve multiple material layers with uncertain parameters we resort to a more robust strategy. In particular, the nonlinear coupling is solved by Newton iterations. We consider the full system of equations which, after discretization in space and time, can be written as 

$$
\mb{F}(\mb{Z},\mb{u}^s, \boldsymbol{\phi}_Q,  \boldsymbol{\phi}_M,  \boldsymbol{\phi}, \mb{S}, \mb{P},  \mb{u}^D )=\mathbf{0}
$$

where $\mb{Z},\mb{u}^s, \boldsymbol{\phi}_Q,  \boldsymbol{\phi}_M,  \boldsymbol{\phi}, \mb{S}, \mb{P}, \mb{u}^D$ are vectors containing the nodal or cell values of the discrete unknowns: the position of the grid nodes $z_i$, the volumetric fraction of precipitated quartz, the mechanical and effective porosity, the sedimentary load, the fluid pressure and the Darcy velocity.

If we set $\mb{X}=\left[\mb{Z},\mb{u}^s, \boldsymbol{\phi}_Q,  \boldsymbol{\phi}_M,  \boldsymbol{\phi}, \mb{S}, \mb{P}, \mb{u}^D\right]$ the $j-th$ iteration of Newton's method reads

\begin{eqnarray*}
&&\mb{J_F(X^{(j)})}\boldsymbol{\delta}\mb{X^{(j)}} = -\mb{F(X^{(j)})}\\
&&\mb{X^{(j+1)}} = \mb{X^{(j)}} + \boldsymbol{\delta}\mb{X^{(j)}}
\end{eqnarray*}
being $\mb{J_F(X)}$ the Jacobian matrix of $\mb{F(X)}$. For each time step the initial guess is given by the solution at the previous time step. The stopping criterion is a simple test on the grid nodes:

$$\max_{i=1,\ldots,N}\dfrac{|Z_i^{(j+1)}-Z_i^{(j)}|}{h_i^{(j)}}\leq \mathrm{tol}$$

\Lc{}{where ${h_i^{(j)}}$ denotes the size of the $i-th$ cell at the $j-th$ iteration of Newton's method.}

A \Lc{more rigorous}{more stringent} test should in principle account for the norm of the vector of the normalized increments. However, due to the different scaling of the variables it is difficult in practice to select an effective tolerance for the whole vector, while a simpler test on the deformation of the domain is robust and representative of the compaction phenomenon.

Note that the position of the \Lc{$N+1$th}{$N+1$st} node is given as a boundary condition by the paleobathimetry reconstruction.

\begin{remark}
In the present implementation the temperature equation is solved at each Newton iteration separately from the others, after the computation of $\mb{X}^{(j+1)}$. This choice is motivated by the fact that temperature is only weakly influenced by the solution of the other equations in the model. However, the nonlinear solver could be easily modified to include the temperature equation. 
\end{remark}


\section{Uncertainty Quantification of  the discontinuous outputs of the compaction model}\label{sec:met}

As we have already mentioned, our goal is to extend the previous works \cite{Colombo16_HJ,feal:compgeo,lever.eal:inversion},
to be able to characterize the uncertainty on the quantities of interest 
computed by the model for multi-layered basins. \Lc{described in the previous section, 
either state variables or quantities obtained post-processing the state variables.}{} 
\La{Throughout this paper, we are only concerned with Uncertainty Quantification 
for quantities at final time, i.e., $T=$today,
but the procedure explained in the following applies verbatim to quantities at any intermediate
time-step $t<T$, under minimal assumptions that we will discuss at the end of this section.}

Specifically, we suppose that the compaction model has $\myn$ uncertain parameters%
\Lc{(we will discuss them in details later on)}{}, which we denote by $p_1,p_2 \ldots,p_{\myn}$, 
\Lc{we collect all of them}{and collect them} in an ${\myn}$-dimensional random vector, $\pp = [p_1,p_2 \ldots,p_{\myn}]$.
We assume that each uncertain parameter $p_n$ can range in 
the interval $\Gamma_n = [a_n,b_n]$, with $a_n,b_n \in \Rset$, and 
that each value is equally likely, i.e, we model each $p_n$ 
as a uniform random variable, with probability density 
function $\rho_n(p_n) = \frac{1}{b_n - a_n}$. 
The random vector $\pp$ will therefore take values in the ${\myn}$-dimensional
hyper-rectangle $\Gamma = \Gamma_1 \times \Gamma_2 \times \ldots \times \Gamma_{\myn}$,
$\pp \in \Gamma \subset \Rset^{\myn}$; 
furthermore, it is reasonable to assume that such random variables are independent,
therefore the joint probability density function of $\pp$ will be 
$\rho(\pp) = \prod_{n=1}^{\myn} \frac{1}{b_n-a_n}$.
Denoting by $f=f(\pp): \Rset^{\myn} \to \Rset$ any output of the compaction model \Lc{(at $T=$today, as already specified)}{at final time $T=$today},
we are interested in computing quantities such as the expected value and the variance of $f$,
\begin{equation}\label{eq:E_and_Var}
\Expct[f] = \int_\Gamma f(\pp) \rho(\pp) d\pp, \quad
\Var[f]=\Expct[f^2] - \Expct[f]^2  
\end{equation}
or its probability density function $\rho_f(\pp)$.

These operations involve multiple evaluations of $f$
for different values of the unknown parameters, and therefore can be very expensive.
Thus, a widely used approach to reduce the computational
cost of this kind of analysis consists in first building a surrogate model for $f$ 
and then obtaining the quantities detailed above by post-processing. \Lc{;
of course, this strategy is effective only if the number of full model evaluations
needed to build the surrogate model is ``small''. To this end, observe that
$f$ is a ${\myn}$-dimensional function, with ${\myn}\gg 1$ typically. Thus, ``naive'' 
approaches that require evaluating $f$ over a cartesian grid of points of $\Gamma$
would be unpractical, because the number of full model evaluations required
by such strategies would grow exponentially with ${\myn}$ (``curse of dimensionality'').
Sophisticated model reduction techniques must therefore be used to circumvent this 
challenge. Here, we follow the previous works}{Following the previous works \cite{feal:compgeo,lever.eal:inversion}}  
\Lc{and consider a}{we consider a sophisticated model reduction technique based on a} global polynomial
surrogate model over the set of parameters, built by a sparse grid approximation formula 
\cite{babuska.nobile.eal:stochastic2,barthelmann.novak.ritter:high,b.griebel:acta,nobile.tempone.eal:aniso,xiu.hesthaven:high}.
However, building an effective sparse grid approximation of state variables such as porosity in the multi-layered
case is not straightforward, due to the fact that $f$ might actually be discontinuous with respect to $\pp$, 
and a dedicated procedure needs to be introduced. 
In the rest of this section, we briefly recall the stardard sparse grids construction,
and then we detail the procedure that we propose to deal with discontinuous response functions. 
In what follows, let $\Nset_+$ and $\Rset_+$ denote the set of strictly greater than 0 integer and real numbers.
Moreover, we will denote by $\eee_k$ the $k$-th canonical vector of $\Rset^{\myn}$,
i.e. a vector with $e_k = 1$ and $e_j=0$ for $j \neq k$ for $k=1,\ldots,{\myn}$, and
by bold numbers, such as $\zzero$ and $\oone$, the vectors of repeated entries
$\zzero = [0,\,0,\, \ldots,0] \in \Nset^{\myn}, \oone = [1,\,1,\,\ldots,1] \in \Nset^{\myn}$.

\subsection{Sparse grid approximation}\label{sec:sparse-grids-basics}

The basic block to construct a sparse grid approximation is a sequence of univariate
interpolant polynomials over each $\Gamma_n$ (in our case, Lagrange polynomials, as detailed below), 
indexed by $j \in \Nset$, which thus denotes the univariate interpolation level. 
For any direction $n=1,\ldots,{\myn}$ of the parameter space, consider a set of $m(j)$ points, 
$\mcP_{m(j)} =\{ \param_{j,1},\param_{j,2},\ldots,\param_{j,m(j)} \} \subset \Gamma_n$,
where $m(\cdot): \Nset \rightarrow \Nset$ is a function
called ``level-to-nodes'' function that specifies the number of points
to be used at each interpolation level, such that $m(0)=0,m(1)=1$ and $m(j)\geq (j-1)$.
\Lc{the sequence of sets $\mcP_{m(1)},\mcP_{m(2)},\ldots$ is said to be ``nested''
if the set of points at level $j$ is a subset of the set of points at level $j+1$,
$\mcP_{m(j)} \subset \mcP_{m(j+1)}$, and non-nested otherwise.}{}
The points should be chosen according to the probability density functions $\rho_n$, 
see e.g. \cite{babuska.nobile.eal:stochastic2}; Gauss--Legendre \Lc{(non-nested)}{} and 
Chebyshev/Clenshaw--Curtis \Lc{(nested)}{} points are suitable
choices for uniform random variables\Lc{, while equispaced points should not be chosen
since they are well-known to suffer of accuracy issues (e.g. the Runge phenomenon), }{}
cf. e.g. \cite{trefethen:comparison,trefethen:book2013}.
\Lc{Observe that the family of points could in principle be different on each direction $n$,
i.e., one might want to use Gauss--Legendre points for a specific random parameter and
Clenshaw--Curtis points for another one, but in this work we will make the 
same choice for all the parameters and therefore we write $\param_{j,k}$ instead of
$\param_{j,k,n}$ and $\mcP_{m(j)}$ instead of $\mcP_{m(j)}^n$.}{}

For a fixed level $j$ and set of points $\mcP_{m(j)}$, 
we denote by $\mcU_n^{m(j)}[g]$ the standard Lagrangian interpolant operator of 
a continuous function $g \in C(\Gamma_n)$, and furthermore we define the detail operator $\Delta_{n}^{m(j)}$ as
\begin{equation}\label{eq:delta_def}
\Delta_{n}^{m(j)}[g](\param) = \mcU_n^{m(j)}[g](\param) - \mcU_n^{m(j-1)}[g](\param).  
\end{equation}
Observe that building the interpolant (and hence the detail operator) requires evaluating the function $g$
at the points of $\mcP_{m(j)}$. We also set $\mcU_n^{0}[g](\param) = 0$.

Next, we introduce the multivariate counterparts of $\mcP_{m(j)}, \mcU^{m(j)}_n$ and $\Delta^{m(j)}_n$.
For any $\jj \in \Nset^{\myn}$,  we write $m(\jj)$ to signify (with a slight abuse of notation)
the vector $[m(j_1)\,m(j_2)\,\ldots m(j_{\myn})]$. 
We begin by letting $\mcG_{m(\jj)}$ denote the cartesian grid with 
$m(j_1) \times m(j_2) \times \ldots \times m(j_{\myn})$ points,
$
\mcG_{m(\jj)} = \mcP_{m(j_1)}^1 \times \mcP_{m(j_2)}^2 \times \ldots \times \mcP_{m(j_{\myn})}^{\myn}.
$ 
Upon evaluating any given multivariate continuous function $g:\Gamma \rightarrow \Rset$
over the points of such grid, we obtain the multivariate tensor Lagrange interpolant of $g$, 
\begin{equation}\label{eq:tensor_interp_def}
  \mcT^{m(\jj)}[g](\pparam) = ( \mcU_1^{m(j_1)} \otimes \ldots \otimes \mcU_{\myn}^{m(j_{\myn})} )[g](\pparam). 
\end{equation}
Note that $\mcT^{m(\jj)}[g]=0$ whenever a component of $m(\jj)$ is zero.
Finally, we introduce the so-called hierarchical surplus operators by taking tensor product of details operators,
\begin{equation}\label{eq:ddelta_def}
  \DDelta^{m(\jj)}[g](\pparam) 
  = ( \Delta_1^{m(j_1)} \otimes \ldots \otimes \Delta_{\myn}^{m(j_{\myn})})[g](\pparam).
\end{equation}   
The sparse grid approximation of $f$ is then defined as a sum of hierarchical surpluses. Namely, we consider
a downward closed set%
\footnote{
A downward closed set (also referred to as ``lower set''), is a set $\mcI \subset \Nset_+^{\myn}$ such that
\[
\forall \iii \in \mcI, \quad \iii - \eee_k \in \mcI \text{ for } k \in \{1,2,\ldots,{\myn}\} \text{ such that } i_k>1.
\]
} $\mcI \subset \Nset_+^{\myn}$ and define 
\begin{equation}\label{eq:sparse_grids_def}
\mcS_{\mcI}[f](\pparam) \,\, = \,\, \sum_{\jj \in \mcI} \DDelta^{m(\jj)}[f](\pparam)  
\,\, = \,\, \sum_{\jj \in \mcI} c_{\jj} \mcT^{m(\jj)}[f](\pparam),
\quad 
c_{\jj} = \mathop{\sum_{\kk \in \{0,1\}^{\myn}}}_{(\jj+\kk) \in \mcI}(-1)^{|\kk|},  
\end{equation}
where the second equality is known as ``combination technique form'' of the sparse grid and can be obtained by combining
equations \eqref{eq:delta_def}, \eqref{eq:ddelta_def} and \eqref{eq:tensor_interp_def} 
with the first expression in \eqref{eq:sparse_grids_def}, see e.g. \cite{wasi.wozniak:cost.bounds}.
The ``combination technique form'' can be useful in practical implementations. Observe also that many $c_{\jj}$ will be zero:
more precisely, if $\mcI$ is downward closed then $c_\jj$ is zero whenever $\jj+\oone \in \mcI$.

The set $\mcI$ in \eqref{eq:sparse_grids_def} prescribes the hierarchical
surpluses to be included in the sparse grid, and should be chosen to give good approximation properties 
while keeping the number of points in the sparse grid to a minimum (remember
that one full PDE model solve per point is required to build the sparse grid).
Observe that choosing 
$
\mcI = \left\{ \jj \in \Nset_+^{\myn} : \max_{n=1\ldots,{\myn}} i_n \leq w \right\} \quad \text{for } w \in \Nset
$
would result in 
$
\mcS_{\mcI}[f](\pparam) = \mcT^{m(\bm{w})}[f](\pparam),
$
i.e., a multivariate tensor approximation of $f$ over a cartesian grid with $m(w)$ points along each
direction $\Gamma_n$, that \Lc{we have already ruled out as}{is} impractical due to the fact that it would require $m(w)^{\myn}$ PDE solves.  
Common choices for $\mcI$ are 
\begin{equation}\label{eq:iso-set}
\mcI = \left\{ \jj \in \Nset_+^{\myn} : \sum_{n=1}^{\myn} (i_n-1) \leq w \right\} \quad \text{for } w \in \Nset
\end{equation}
and the anistropic counterpart
\begin{equation}\label{eq:aniso-set}
\mcI_{\aalpha} = \left\{ \jj \in \Nset_+^{\myn} : \sum_{n=1}^{\myn} \alpha_n(i_n-1) \leq w \right\} \quad \text{for } w \in \Nset  
\end{equation}
where $\aalpha = [\alpha_1,\ldots,\alpha_{\myn}] \in \Rset^{\myn}$, $\alpha_n > 0$, are coefficients that are chosen so that a higher
polynomial degree is allowed along the directions deemed more important
(the more important the random parameter, the smaller the corresponding coefficient $\alpha_n$), 
see e.g. \cite{feal:compgeo,nobile.tempone.eal:aniso,back.nobile.eal:comparison}.
More advanced strategies for the selection of the set $\mcI$ are available in literature:
for instance, in \cite{nobile.eal:optimal-sparse-grids,back.nobile.eal:optimal,back.nobile.eal:lognormal} 
an algorithm for optimal sparse grids construction for elliptic PDEs with random coefficients is discussed,
while ad adaptive algorithm based on ``a-posteriori'' refinement has been discussed in
\cite{gerstner.griebel:adaptive,schillings.schwab:inverse,chkifa:adaptive-interp,nobile.eal:adaptive-lognormal,narayan:Leja}.
In this work we only discuss results obtained with anisotropic sets of the form \eqref{eq:aniso-set}
(the specific choice of $\aalpha$ will be detailed later, in Section \ref{sec:synt-case}). 
We have also performed some tests with adaptive construction of $\mcI$ and obtained
analogous results (not shown).

We conclude this introductory \Lc{sections}{section} by briefly recalling how relevant Uncertainty Quantification
information can be extracted from a sparse grid approximation by suitable post-processing.
First, let $\mcH^{m}_{\mcI}$ be the set of collocation points used to build the sparse grid approximation
$\mcS_{\mcI}[f]$; 
one can then associate a sparse grid quadrature formula $\mcQ_{\mcI(w)}^m[\cdot]$ to 
$\mcS_{\mcI}[f]$, i.e. 
\begin{equation}\label{eq:sparse-quad}  
\forall \,\, f \in C^0(\Gamma), \,\,
\int_\Gamma f(\pp) \rho(\pp)d\pp  \approx
\int_\Gamma \mcS_{\mcI(w)}[f](\pp)\rho(\yy)d\yy =
\mcQ_{\mcI(w)}[f] =  \sum_{\pp_j \in \mcH^{m}_{\mcI}} f(\pp_j)\varpi_j,
\end{equation}
for suitable weights $\varpi_j \in \Rset$. 
As a consequence, upon building the sparse grid approximation of $f$, one can immediately compute an approximation
of its expected value and higher moments introduced in \eqref{eq:E_and_Var}.

Second, by some algebraic manipulations of the combination technique form of the sparse grid approximation
$\mcS_{\mcI}[f](\pparam)$, cf. eq. \eqref{eq:sparse_grids_def}, it is possible to compute 
the so-called Sobol indices of $f$, which quantify which percentage of the total variability 
of $f$ is due to the variability of each random parameter and to their products, 
in the spirit of a variance decomposition analysis see e.g. \cite{Saltelli2006,lemaitre:book,sobol93,sudret:sobol}.
We refer the interested reader to \cite{feal:compgeo,lorenzo:thesis} for more details on how to 
obtain the Sobol indices of $f$ starting from its sparse grids approximation $\mcS_{\mcI}[f](\pparam)$;
see also \cite{constant.eldred.phip:pseudospec}.

\subsection{Approximation of discontinuous outputs\Lc{: a two-steps surrogate model}{}}\label{sec:disc-output}


\Lc{As already discussed in the introduction, in this work we focus on
sedimentary basins composed by layers of different geomaterials (or
lithological units). While some of the variables are continuous across
material interfaces (e.g., temperature and pressure vertical
distributions) other key quantities (e.g., porosity and permeability)
may significantly differ as a function of the geomaterial. At the same
time, the depth of the interface between lithological units is random,
as it typically depends on the values of the random parameters.}{In this work we focus on
sedimentary basins composed by layers of different geomaterials whose interfaces position is random,
as it typically depends on the values of the random parameters.} As a
consequence different kind of geomaterials may be found at a given
depth $z^*$ at time $t$ for different realizations of the random
parameters. This phenomenon is exemplified in Figure
\ref{fig:smooth-and-non-smooth}-left, which shows two vertical
distributions of porosity obtained from the direct solution of the
mathematical model described in Section \ref{sec:model} 
for two different realizations of the random parameters.
The jumps in the porosity profiles are symptoms of geomaterial changes, and 
it is clearly visible from the figure that interfaces occur at different dephts depending on the random parameter. 
This implies that the porosity (and other quantities of interest as well) 
at a fixed depth $z^*$ will be exhibit a discontinuous behavior with respect to 
the random parameters $\pp$.

\Lc{As a consequence, the quality of a standard sparse grids approximation for $f(\pp)$
will be dramatically affected, since sparse grids interpolants are effective only if $f(\pp)$ is a smooth function of the parameters $\pp$. The Uncertainty Quantification procedure for single-material
basins developed in the previous work \cite{feal:compgeo} cannot be applied
verbatim in this context, or, more precisely, it can be applied only at those depths
for which the same material occurs regardless of the value of the random parameters.
Similarly, with a view to inverse basin modeling, the procedure we proposed in \cite{lever.eal:inversion}
based on \cite{feal:compgeo} could only be applied upon discarding all data collected at depths
for which different materials might occur. This hampers the possibility to apply our methodology in significant real world scenarios, particularly when thin layers appear in the stratigraphic sequence.}

\Lc{In this work we propose an original procedure to circumvent this issue. Our proposed methodology}{As a consequence, the quality of a standard sparse grid approximation for $f(\pp)$ will be dramatically affected, since sparse grids interpolants are most effective when $f(\pp)$ is a smooth function of the parameters $\pp$. To circumvent this issue, we propose a methodology that} stems from the observation that while quantities at a fixed depth will be likely discontinuous with respect to the parameters $\pp$,
the position of the interface between two materials typically depends smoothly on $\pp$. \Lc{We give numerical evidence of this fact}{This can be observed} in Figure \ref{fig:smooth-and-non-smooth}-right, where we display the variation of the position of the material interfaces as a function of two random parameters,
namely the sandstone \Lc{and shale vertical compressibilities}{vertical compressibility and the parameter $k_2$ of the porosity-permeability law}
(refer to Section \ref{sec:synt-case} for more details). This suggests to create a surrogate model
for the state variables with a two-steps approach: first, a sparse-grid approximation
for the position of each interface is computed, and then the state variables are
approximated within each homogeneous lithology.

\begin{figure}[t]
  \centering
  \includegraphics[width=0.44\linewidth]{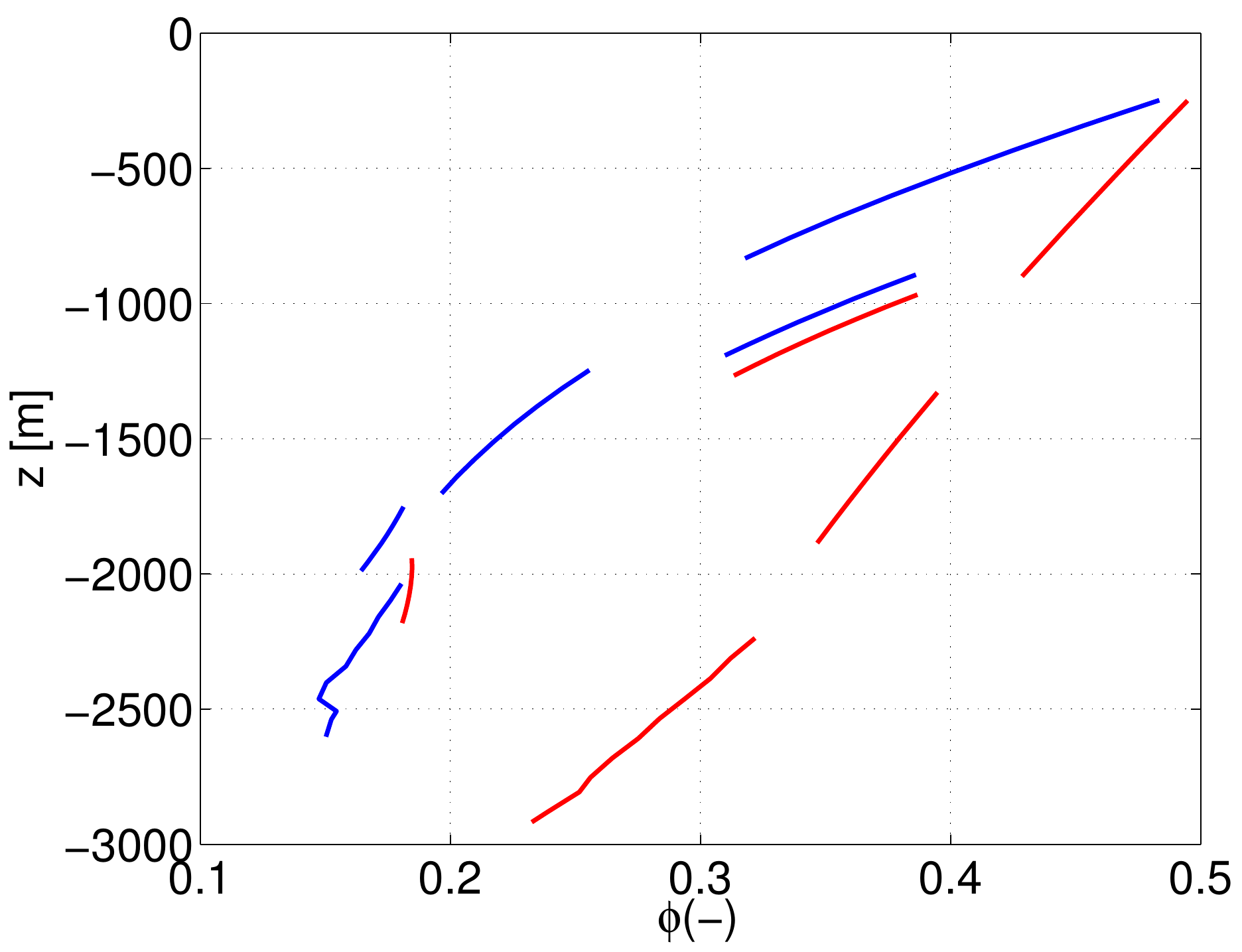}  
  \includegraphics[width=0.54\linewidth]{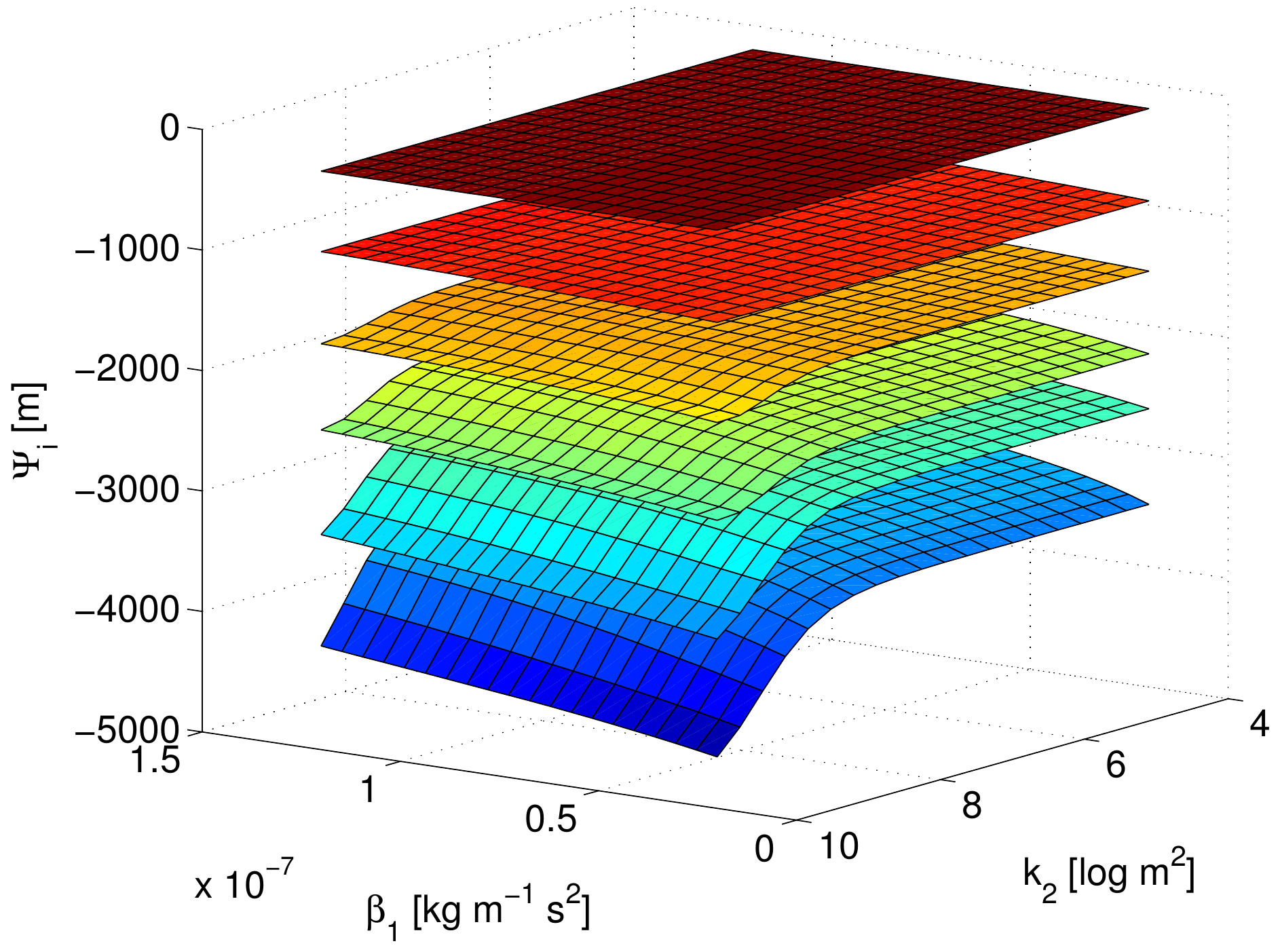}
  \caption{Left: two porosity $z$-profiles (a blue one and a red one) corresponding
    to two different realizations of the
    random parameters $\pp$. The porosity jumps between consecutive
    layers, and the position of the interfaces depends on the
    realizations of $\pp$, leading to discontinuity in the function
    $f(\pp) = \phi(z,\pp)$ for fixed $z$. 
    Right: depth of the interfaces
    as a function of $\pp$ in a case where the layers are composed
    by two alternating materials \Lc{whose mechanical compaction parameters
    (parameters $\beta$ in \eqref{eq:mec-comp})
    are considered random, $\pp=[\beta_1 \, \beta_2] \in \Gamma$.}{
    and $k_2$ and $\beta_1$, cf. eq. \eqref{eq:k-phi} and \eqref{eq:mec-comp}, 
    are considered to be random parameters, $\pp=[\beta_1 \, k_2] \in \Gamma$.} 
    The position of the interfaces varies smoothly over $\Gamma$.}
  \label{fig:smooth-and-non-smooth}
\end{figure}


\subsection{A two-steps surrogate model}\label{sec:disc-output}

\begin{figure}[tb]
  \centering
  \includegraphics[width=0.85\linewidth]{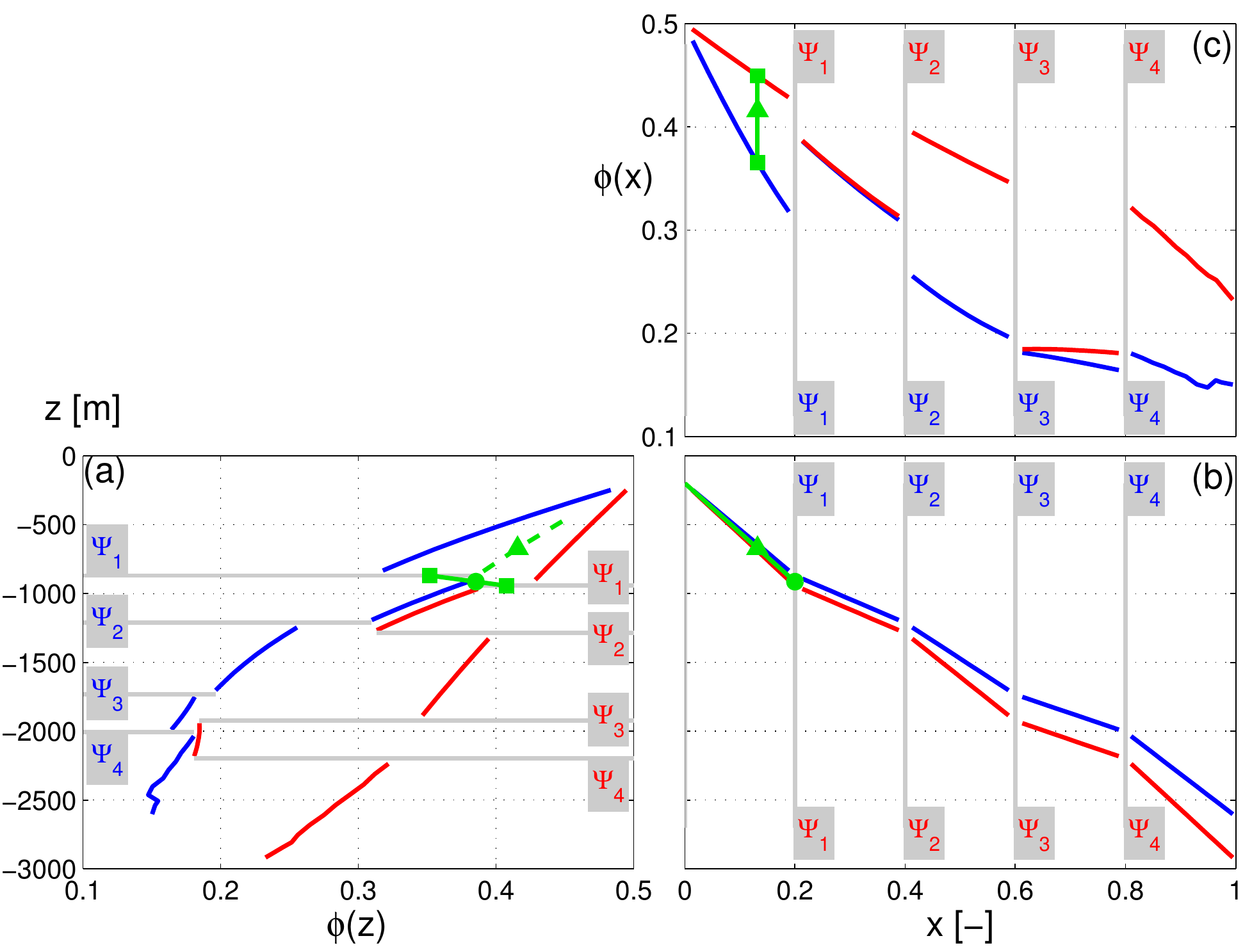}  
  \caption{Cartoon of the two-steps procedure to construct a sparse grid surrogate model of a discontinuous quantity.}
  \label{fig:z-phi-x-plot}
\end{figure}


We detail the procedure referring to the cartoon in Figure \ref{fig:z-phi-x-plot}, which is composed by
three panels, (a),(b), and (c), starting from the bottom-left corner and proceeding counter-clock-wise. 
We focus here on the porosity as an example, but the procedure can be applied verbatim to any 
other discontinuous output of interest.
Figure \ref{fig:z-phi-x-plot}-(a) shows in blue and red two porosity-depth profiles 
at $T$=today, corresponding to two different realizations of the uncertain parameters,
and in green a portion of the porosity profile predicted with our two-steps procedure.

We begin by denoting by $\intpos_k(\pp)$ the depth of the $k$-th interface 
as a function of the random parameters, 
where $k=0$ denotes the seafloor and $k=K$ denotes the bottom of the basin.
Note that our procedure assumes that the number of layers does not
depend on the random parameters.
The seafloor is the upmost interface, i.e., the top of the most recent layer, 
around (fixed to) $-200$m in Figure \ref{fig:z-phi-x-plot}-(a), while the bottom of the basin
is the bottom of the oldest layer, whose depth depends on $\pp$ and is around $-2500$m for 
the blue profile and around $-3000$m for the the red profile. 
As already discussed, the interfaces among geomaterials in the blue and red profiles
occur at different depths. 

The key idea to create the sparse \Lc{grids}{grid} approximations of the state variables 
is to operate on a layer-by-layer basis, introducing a piecewise linear mapping $\hat{x}=F(z,\pp)$, 
that projects the state variables profiles from the physical domain 
into a reference domain where the discontinuities are aligned, or, in other words, where the
interfaces between \Lc{geometerials}{geomaterials} are forced to occur at the same location for every realization of the
random parameters.
Crucially, since the position of the interfaces is random, the map $F$ is random as well, i.e.,
it depends on $\pp$. Figure \ref{fig:z-phi-x-plot} shows two realizations of such map in panel (b) 
and the aligned porosity profiles in Figure \Lc{}{\ref{fig:z-phi-x-plot}}(c).

In details, the piecewise linear mapping is defined so that the interfaces 
are aligned at prescribed values $\hat{x}_k$: 
the seafloor is located at $0 = \hat{x}_0 = F(\intpos_0(\pp),\pp)$, 
the bottom of the basin is located at $1 = \hat{x}_K = F(\intpos_K(\pp),\pp)$,
and intermediate interfaces are located at $\frac{k}{K} = \hat{x}_k = F(\intpos_k(\pp),\pp)$.
Incidentally, note that the reference domain needs not be $[0,1]$ and could be also e.g $[0,K]$, with $\hat{x}_k=k$. 
Assuming for a moment that $\intpos_0(\pp), \intpos_1(\pp),\ldots,\intpos_K(\pp)$ are known for a given $\pp$, then
the map $F$ can be written as
\[
\hat{x} 
	= z \frac{\hat{x}_{k+1}-\hat{x}_k}{\intpos_{k+1}(\pp)-\intpos_k(\pp)} 
	+ \frac{\hat{x}_{k}\intpos_{k+1}(\pp)-\hat{x}_{k+1}\intpos_k(\pp)}{\intpos_{k+1}(\pp)-\intpos_k(\pp)} 
\text{ for } \intpos_k(\pp) \leq z \leq \intpos_{k+1}(\pp).
\]
In the reference domain, the interfaces are aligned, therefore for a fixed $\hat{x}$ the 
state variables will always be related to the same layer: therefore, it is reasonable to assume that
in the reference domain the state variables at fixed $\hat{x}$ depend smoothly on the parameters,
and can be effectively approximated by a sparse \Lc{grids}{grid} approximation. 
In other words, we are implicitely assuming that
the behavior of the state variables at any fixed fraction of depth of each layer depends
smoothly on the parameters.

Clearly, $\intpos_0(\pp), \intpos_1(\pp),\ldots,\intpos_K(\pp)$ are not known a-priori;
however, the functions $\intpos_k(\pp)$ can be assumed to be smooth with respect to $\pp$
as we discussed earlier, therefore they can be effectively approximated by sparse grids.
Thus, the layer-by-layer state variables approximation algorithm can be summarized as follows:
\begin{enumerate}
\item Construct the sparse grid approximation of the position of the interfaces between layers,
  $\mcS^m_{\mcI}[\intpos_k](\pp)$ for every $k=0,1,\ldots,K$;
\item For every couple of values $z^*,\pp^*$ at which an approximation of $\phi(z^*,\pp^*)$ is sought (for instance, the location marked
  by a green triangle in panel (a)):
  \begin{enumerate}
  \item approximate by sparse grids the position of the interfaces for $\pp^*$, $\intpos_k(\pp^*) \approx \mcS^m_{\mcI}[\intpos_k](\pp^*)$:
    in panel (a), the green circle marks the predicted depth at which the first interface will occur, given the 
    location of the first interface for the blue and red profile, marked with squares;
  \item use these values to derive the expression of the piecewise linear mapping $\hat{x}=F(z,\pp^*)$ and compute
    $\hat{x}^* = F(z^*,\pp^*)$: in panel (b), the green circle defines the first portion of the new realization of the piecewise
    linear mapping, that can be used to compute the reference location $\hat{x}^*$ corresponding to the physical location $z^*$ 
    of the green triangle in panel (a);
  \item compute the sparse grid approximation of the porosity at $\hat{x}^*$,
    $\phi(\hat{x}^*,\pp^*) \approx \mcS^m_{\mcI}[\phi(\hat{x}^*,\cdot)](\pp^*)$.
    In panel (c), a green triangle marks the predicted value of porosity at $\hat{x}^*$ given the 
    porosity values at the same location for the blue and red profile, marked with squares.
    Observe that some kind of finite element interpolation for the profiles computed with
    the full model (as the one detailed in Section \ref{section:discretization}) is needed 
    as a preliminary step to compute at $\hat{x}^*$ the values of porosity  
    that will be used to build $\mcS^m_{\mcI}[\phi(\hat{x}^*,\cdot)](\pp^*)$;
  \item approximate $\phi(z^*,\pp^*) \approx \mcS^m_{\mcI}[\phi(\hat{x}^*,\cdot)](\pp^*)$.
  \end{enumerate}
\end{enumerate}

\subsubsection{Discussion on the surrogate model}\label{sec:disc-output}

\La{Before moving to the section showcasing the numerical result, it is worth taking a moment
to highlight and discuss features, extensions and possible limitations of the methodology we
outlined above.}


\La{First, observe that similar approaches to our methodology can be found in 
  \cite{alexanderian1,alexanderian2,liao.zhang:transform}.
  All these works have discussed the advantage that can be obtained by model
  reduction applied to some ad-hoc transformations of the original
  output variables designed to improve the accuracy.}
\Lc{The methodology we just described bears some similarities with the one discussed in the very recent
work \cite{liao.zhang:transform}. 
There, the authors consider a time-dependent advection-diffusion problem, and show that while it is hard to approximate
with sparse grids the dependence on the random parameters of the concentration of a solute at any space-time 
location due to the presence of sharp gradients, the arrival time of a given value of concentration at
the same location is instead more amenable to approximation with sparse grids.
Therefore, the authors suggest a three-step procedure: 
a) pre-process the output of the full-model, i.e., the concentration profiles, 
to obtain the arrival times at the desired spatial locations; 
b) approximate the arrival times for new values of the random parameters with a sparse grid;
c) convert back to obtain the concentration profile at the desired space-time locations. 
Our procedure is different in that we are able to treat discontinuous outputs,
and our first step does not require a pre-process of the porosity profile but is rather
based on a set of independent quantities (the location of the interfaces between layers)
that are interpolated already over the parameter space.}
{Our procedure is different in that it is specifically tuned to tackle problems with discontinuities.}


\La{Concerning possible extensions, as we have already mentioned, 
the procedure that we have detailed can be applied verbatim to any time up to $T=$today.
The only required assumption is that the depositional time of each material be not a random quantity. Indeed, 
if that was the case, the location of the interfaces would be discontinuous with respect to
these random depositional times as well, and this issue would require to be treated analogously to what we propose here for the 
porosity discontinuities.}

\La{We also highlight that the uncertainty quantification methodology proposed 
here can be extended to the outputs of three-dimensional or two-dimensional 
models of compaction processes, by applying it pointwise at all locations where uncertainty quantification may be of interest.
This can be done under the assumption that the same number of layers is encountered along the vertical direction at all locations. 
The procedure cannot be applied in its current implementation in the presence of local variations of 
layering structure which may arise in three-dimensional geological bodies, e.g. in the presence of faults. 
While we envision dedicated remedies that may be formulated for those specific cases, their discussion is beyond the scope of the present contribution.}

\La{Finally, we remark that 
  the choice of the mapping $F$ is arbitrary as long as it is globally invertible and smooth within each layer,
  and thus we have chosen the linear one for convenience.
  Moreover, our approach rests on the assumption that the limiting factor in the accuracy of the surrogate model 
  is actually the fact that the quantity to be interpolated is discontinuous between interfaces, hence its
  composition with a piecewise linear mapping does not reduce the overall smoothness.
}

  



\section{Numerical Results}\label{application}
 In this section we illustrate the application of the approach discussed in section \ref{sec:met} to a couple of synthetic test cases. 
\La{All results reported hereafter have been obtained with Matlab$^{\textregistered}$ 8.5, 
using the sparse grids functionalities provided by \cite{lorenzo:sparse-grids-code}.}


\subsection{Assessment of the numerical solver for the full model}\label{anna}

As we explained in section \ref{section:discretization} we have chosen, for the solution of the coupled nonlinear system of equations, a monolithic approach with Newton iterations. To prove that this approach is more robust than the iterative splitting proposed in previous works we consider a simple test case and compare the solutions obtained with the strategy proposed in \cite{feal:compgeo} and the Newton method. We consider a single layer that is subject to compaction due to its own weight and a constant overburden $S_0$. This simple setup avoids inconsistencies due to different initialization of the new cells added to the domain in the two approaches. 

The main material and geometric properties are summarized in \TAB \ref{table:casepf}.

\begin{table}[htbp]
\centering
 \begin{tabular} {|l|c|l|c|}
 \hline
 Initial thickness & 500$m$ & Initial porosity $\phi_0$ &0.55\\ \hline
 Sea depth $h_{sea} $  & 500$m$ & Compressibility $\beta$  & 4.0e-7\\ \hline
 Rock density $\rho_s$ & 2648$kg/m^3$ & Time span & 2.5$y$\\ \hline
 \end{tabular}
 \caption{Main physical and geometrical parameters for test case \ref{anna}.}\label{table:casepf}
\end{table}

At the initial time $t=0$ the porosity is uniform and the effective stress is set to zero: then, compaction occurs quickly until the equilibrium configuration is reached. \Lc{In this test geochemical reactions are deactivated. Even if a rigorous proof is out of the scope of this work, we want to show by numerical experiments that the iterative splitting, implemented as in \cite{feal:compgeo}, requires more and more iterations as the rock permeability decreases and eventually fails to converge. To this aim we}{We} consider the permeability law (\ref{eq:k-phi}) and the following values for $k_1$ and $k_2$:

$$k_1 =14.9 \alpha + 1.94(1-\alpha),\qquad k_2= 7.7 \alpha + 8 (1-\alpha),  $$

with $\alpha\in[0,1]$ to mimic a variable blend of sandstone and shale. For $\alpha=1$  the sediments are quite permeable and the layer compacts without remarkable overpressure, see figure \ref{pic:aa1}.

\begin{figure}[tb]
 \includegraphics[width=0.33\textwidth]{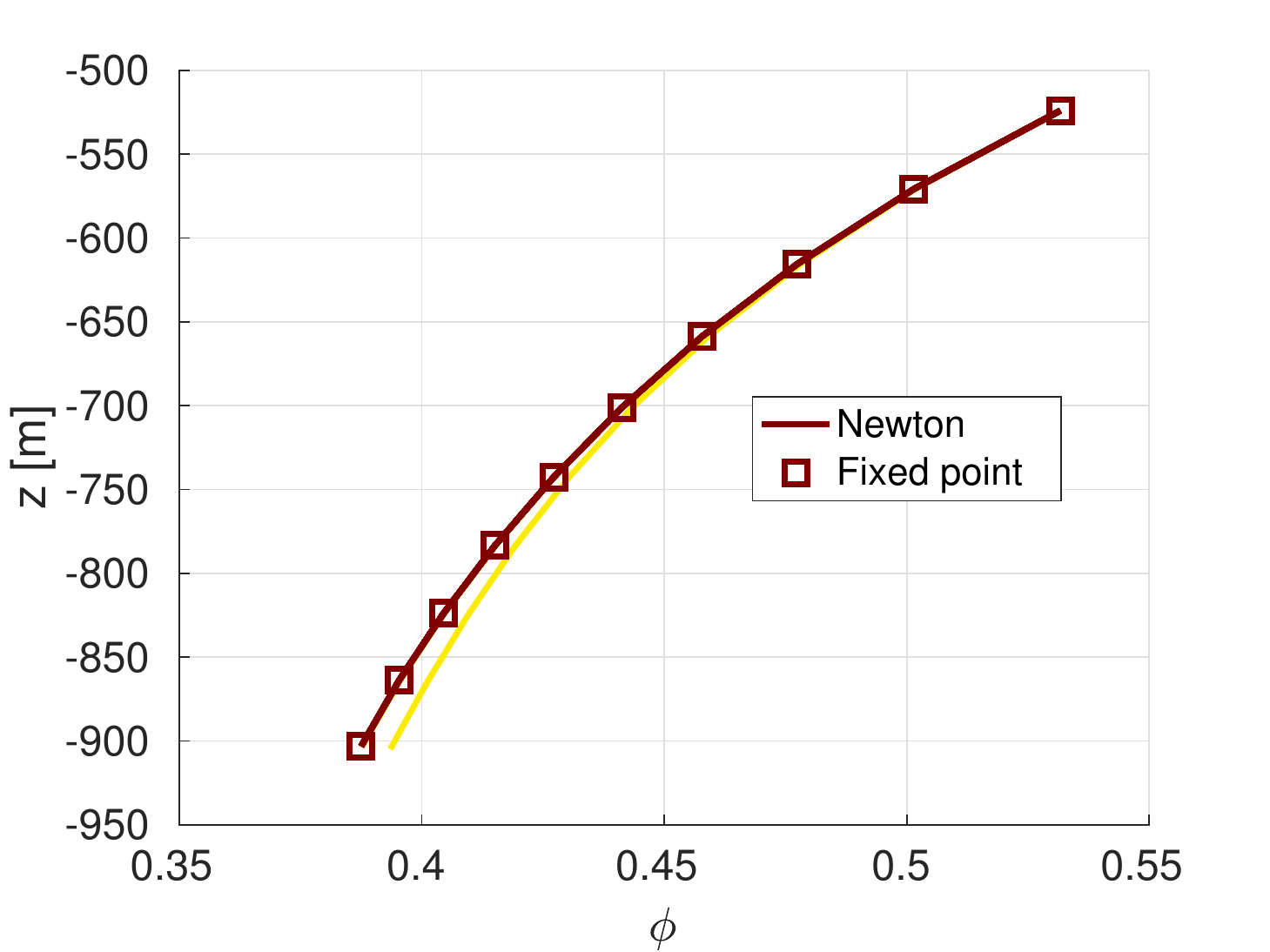}
 \includegraphics[width=0.33\textwidth]{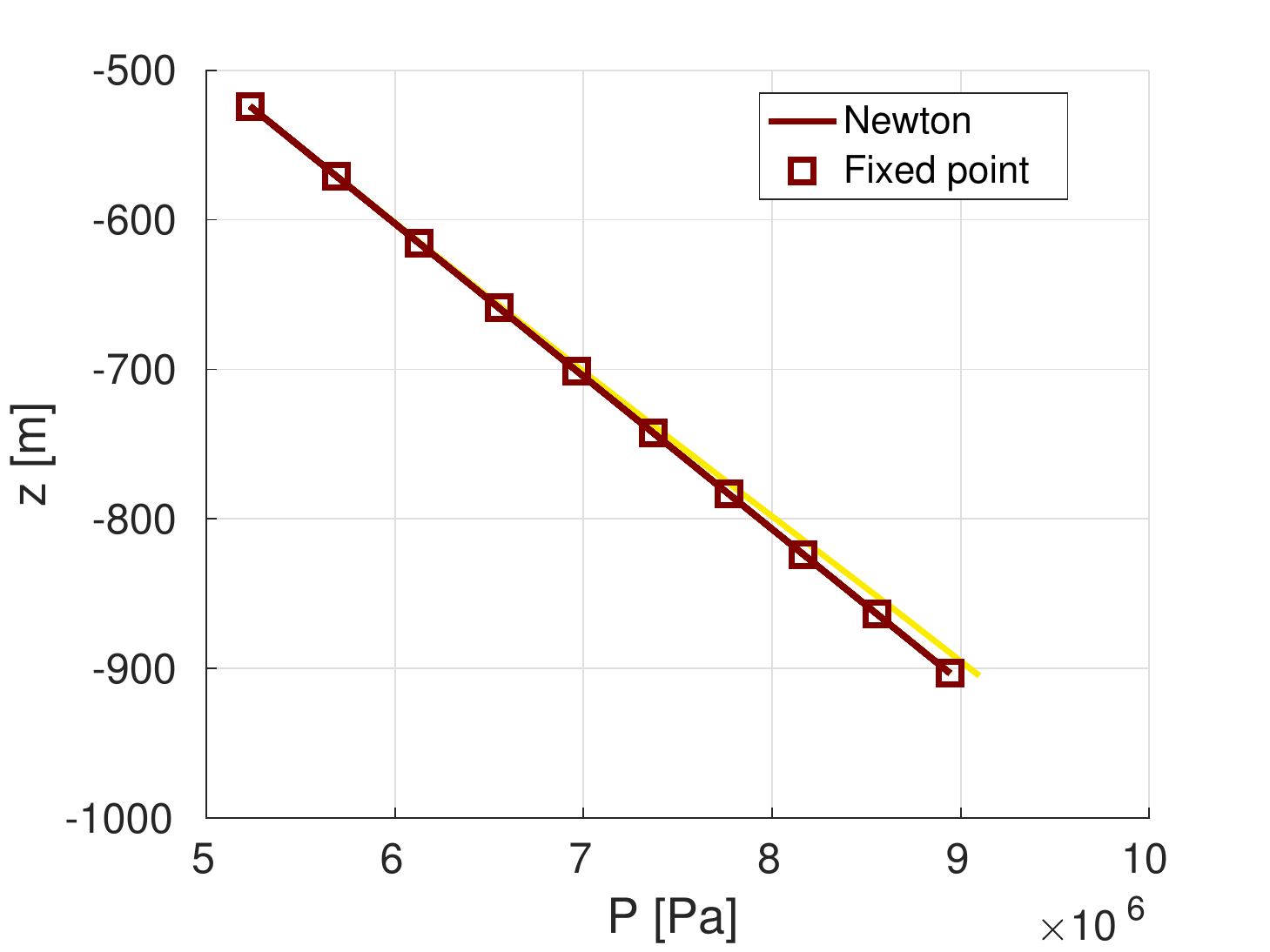}
 \includegraphics[width=0.33\textwidth]{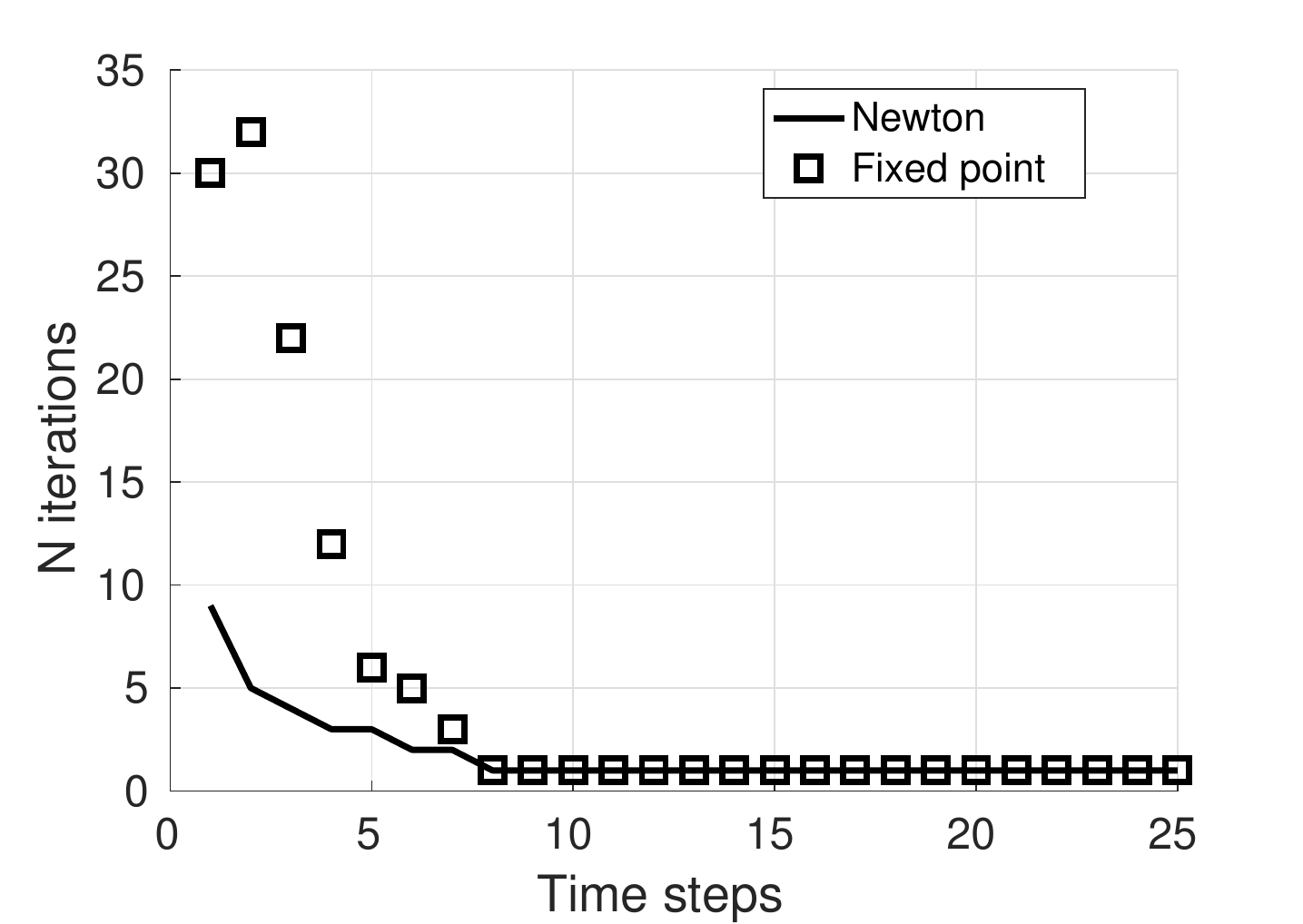}
 \caption{Porosity, pressure and number of iterations for the fixed point and the Newton method for $\alpha=1$. In the first two plots the Newton solution at different times is represented by lines of colors ranging from yellow to red, together with the fixed point solution at the final time.}\label{pic:aa1}
\end{figure}

 As $\alpha$ decreases the behavior changes and we observe overpressure counteracting compaction (figure \ref{pic:aa0p7}). This more complex behavior reflects in a higher number of iterations in particular for the iterative splitting (fixed point) method. Finally for $\alpha=0.5$ the fixed point method fails to converge while, with the same tolerance and stopping criterion we obtain a solution in 8 iterations with the Newton method. This is due to the complex pressure-porosity interaction, shown in figure \ref{pic:aa0p5}: in this \Lc{cases}{case}, due to the low permeability the iterative splitting of \cite{feal:compgeo} is unstable, while we can always converge to the solution with a fully coupled approach.

\begin{figure}[tb]
 \includegraphics[width=0.33\textwidth]{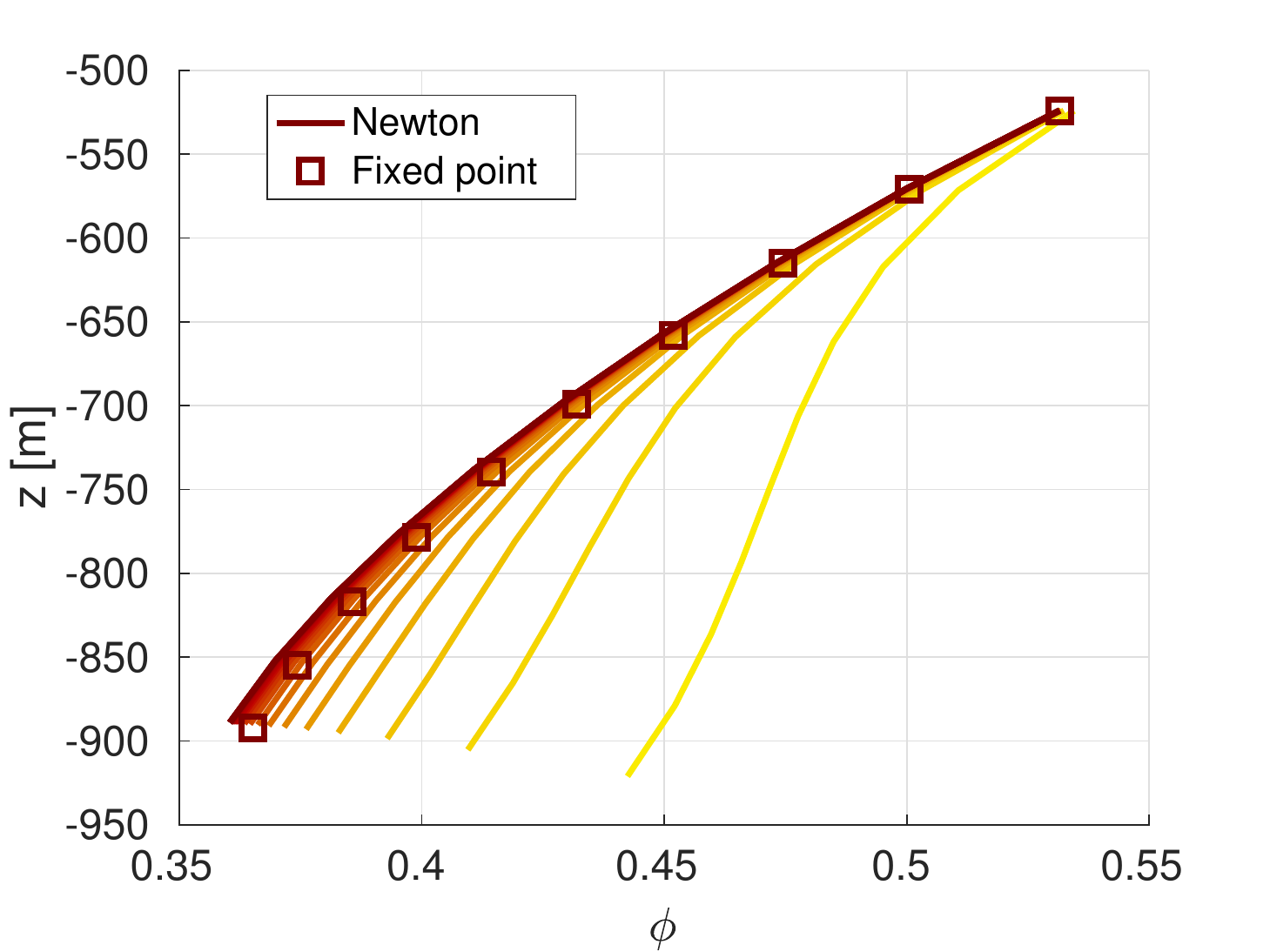}
 \includegraphics[width=0.33\textwidth]{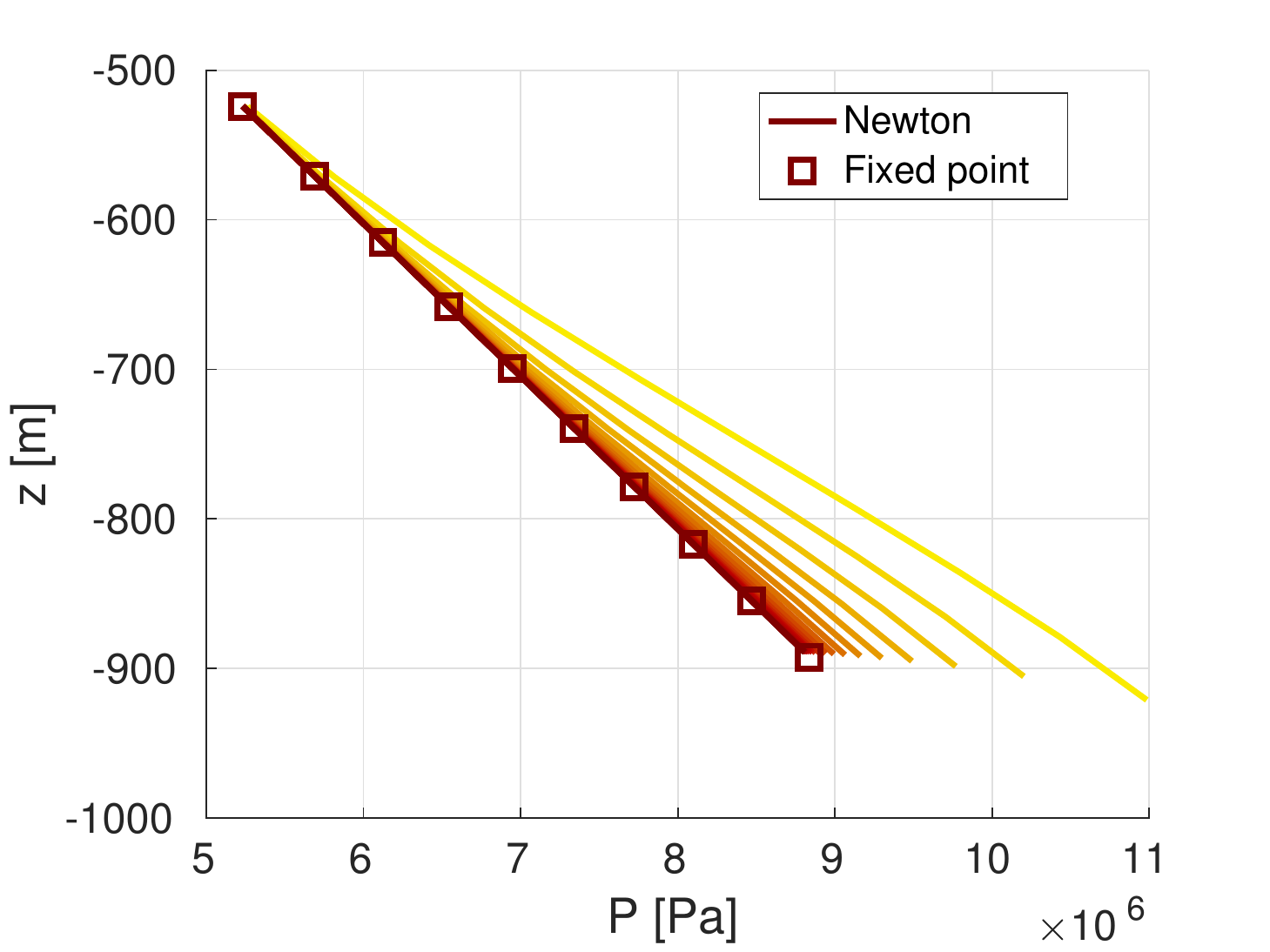}
 \includegraphics[width=0.33\textwidth]{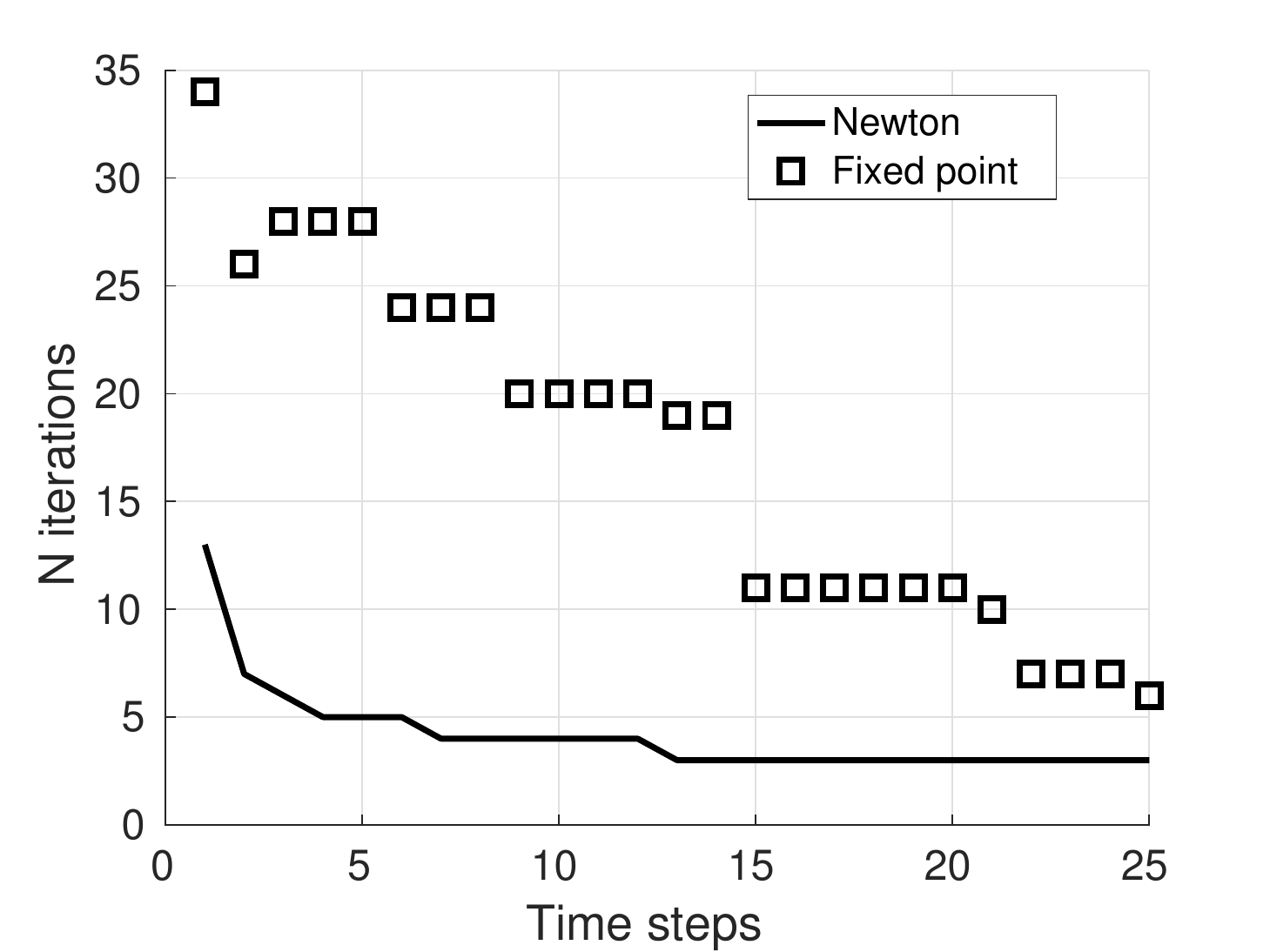}
 \caption{Porosity, pressure and number of iterations for the fixed point and the Newton method for $\alpha=0.7$. In the first two plots the Newton solution at different times is represented by lines of colors ranging from yellow to red, together with the fixed point solution at the final time.}\label{pic:aa0p7}
\end{figure}

\begin{figure}[tb]
 \includegraphics[width=0.33\textwidth]{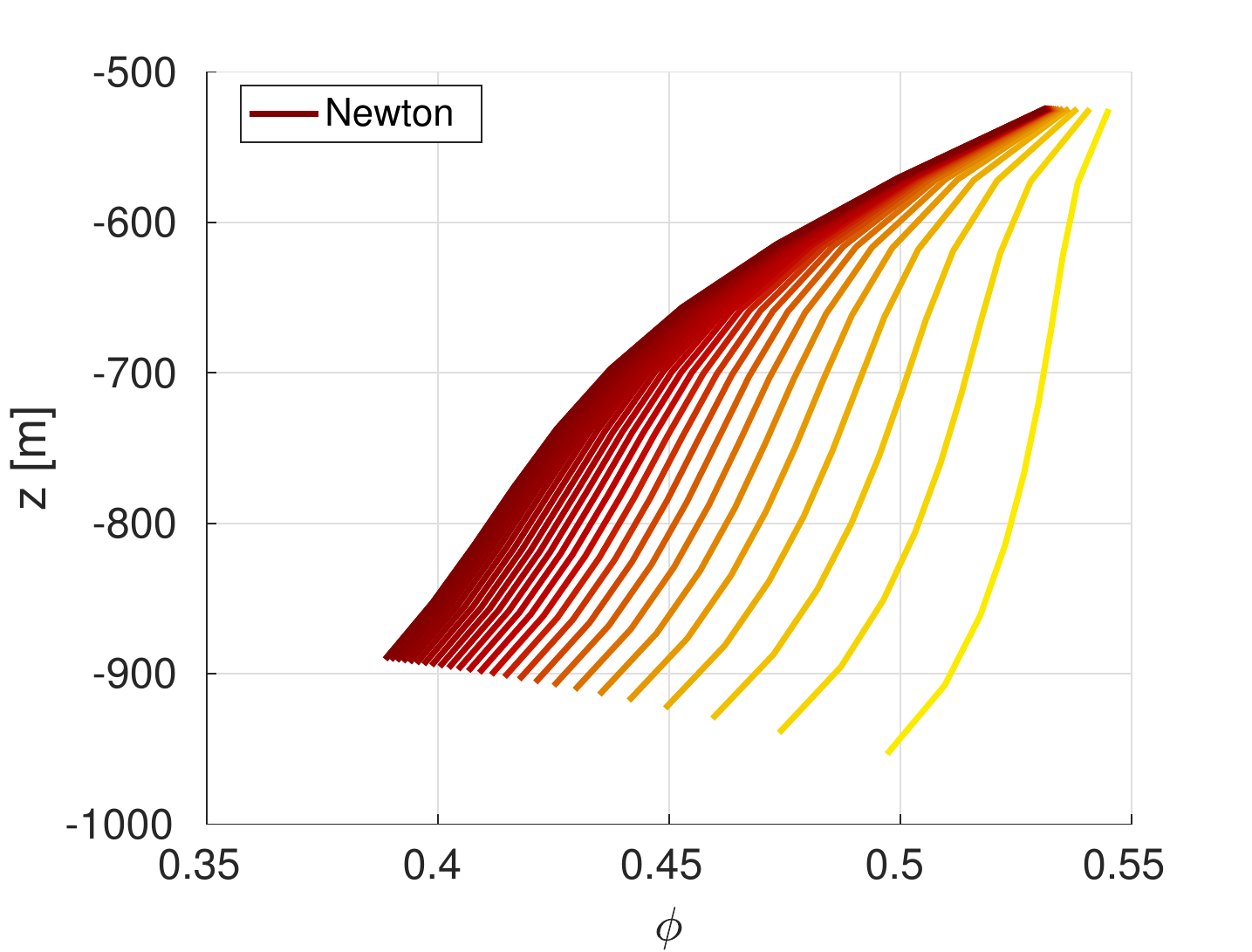}
 \includegraphics[width=0.33\textwidth]{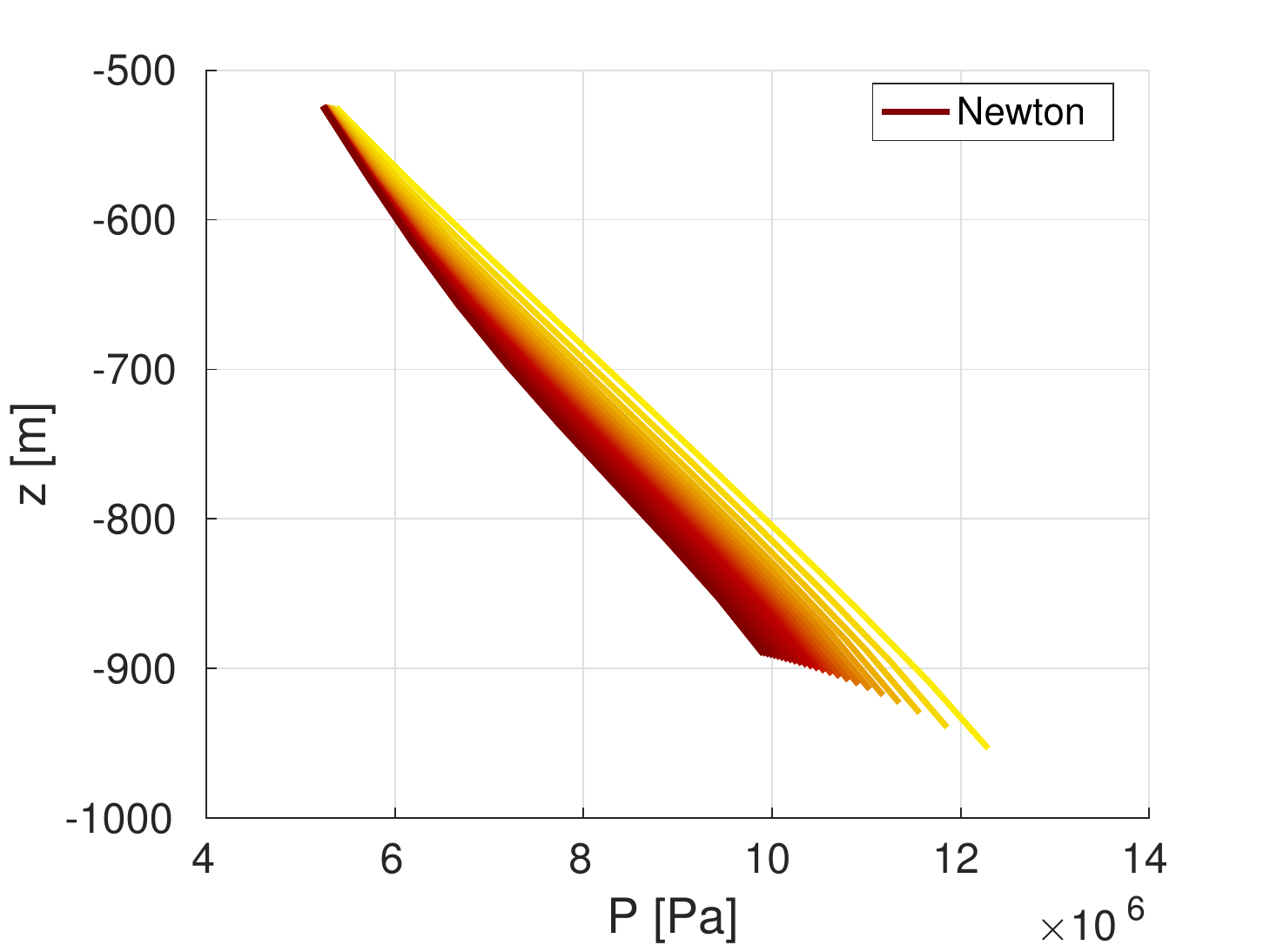}
 \includegraphics[width=0.33\textwidth]{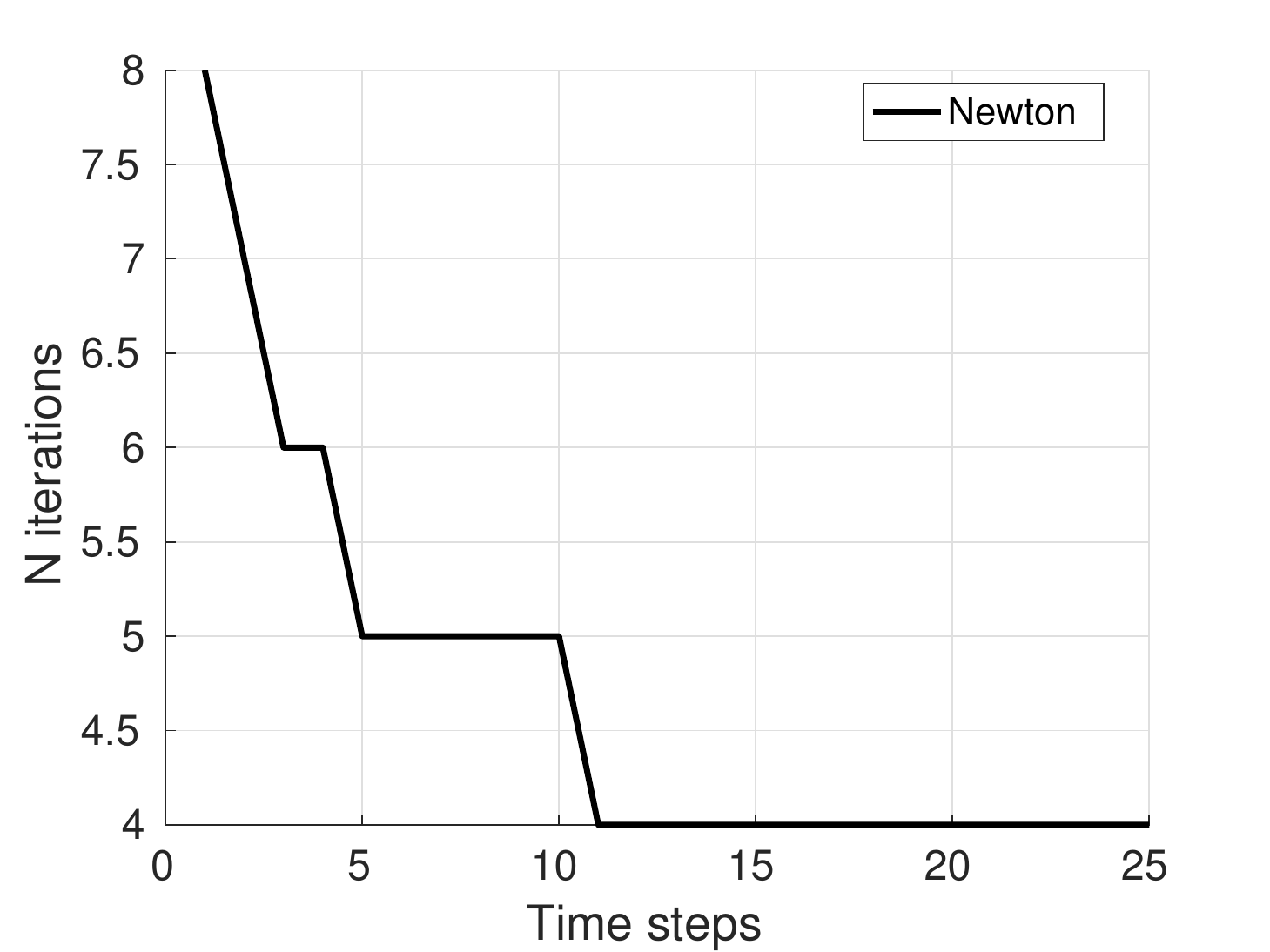}
 \caption{Porosity, pressure and number of iterations for the Newton method for $\alpha=0.5$. In the first two plots the Newton solution at different times is represented by lines of colors ranging from yellow to red.}\label{pic:aa0p5}
\end{figure}

\subsection{Assessment of the Uncertainty Quantification methodology}\label{sec:synt-case}

\begin{table}[tbp]
\begin{center}
\begin{tabular}{|l|c|c|}
\hline
\head{physical meaning} & \head{symbol} & \head{parameter space} \\
\hline
param. for compaction in sandstone & $\beta_{sd}$ & $2\times10^{-8} -  12\times10^{-8}$ \\
\hline
param. for compaction in shale & $\beta_{sh}$ & $2\times10^{-8} -  12\times10^{-8}$ \\
\hline
param. of K-$\phi$ law & $k_2^{sh}$ & 4 - 10 \\
\hline
\end{tabular}
\end{center}
\caption{List of uncertain parameters and associated ranges of variability}
\label{table:uncertain-parameters}
\end{table}

\begin{table}[tbp]
\begin{center}
\begin{tabular}{|p{7.5cm}|c|c|}
\hline
\multicolumn{3}{|c|}{\textbf{SEDIMENTATION CHARACTERISTICS}} \\ 
\hline
\head{parameter} & \head{value} & \head{units} \\
\hline
total sedimentation time & 100 & Ma \\
\hline
constant sedimentation rate & 40 & m $\mbox{Ma}^{-1}$ \\
\hline
depositional order (and time) & 
\hspace{-20pt}
\begin{minipage}{0.15\linewidth}
\begin{enumerate}
\item sand (20)
\item shale (20) 
\item sand (20) 
\item shale (20) 
\item sand (20)  
\end{enumerate}  
\end{minipage}
&  (depositional time [Ma])\\
\hline
\multicolumn{3}{|c|}{\textbf{INITIAL AND BOUNDARY CONDITIONS}} \\ 
\hline
\head{parameter} & \head{value} & \head{units} \\
\hline
bathymetric high ( $h_{sea} $ ) & 200 & m \\
\hline
temperature at basin top ( $T_{top} $ ) & 295.15 & K \\
\hline
heat flux at basin basement ( $G_{T}$ ) & 0.024 & K $\mbox{m}^{-1}$ \\
\hline
Darcy flux at basin basement ( $u^{D}$ ) & 0 & m $\mbox{s}^{-1}$ \\
\hline
\multicolumn{3}{|c|}{\textbf{FLUID CHARACTERISTICS}} \\  
\hline
\head{parameter} & \head{value} & \head{units} \\
\hline
fresh water density  ( $\rho^{l}$ ) & 999 & kg $\mbox{m}^{-3}$  \\
\hline
seawater density  ( $\rho_{sea}^{l}$ ) & 1025 & kg $\mbox{m}^{-3}$  \\
\hline
viscosity ( $\mu^{l}$ ) & $1.001 \times 10^{-3}$ & kg $\mbox{m}^{-1}$ $\mbox{s}^{-1}$   \\
\hline
specific thermal capacity ( $c^{l}$ ) & 4186 & J $\mbox{K}^{-1}$ $\mbox{kg}^{-1}$ \\
\hline
specific heat conductivity ( $\lambda_l$ ) & 0.6 & W $\mbox{K}^{-1}$ $\mbox{m}^{-1}$ \\
\hline
\multicolumn{3}{|c|}{\textbf{POROUS MEDIUM CHARACTERISTICS}} \\  
\hline
\head{parameter} & \head{value} 			 											 & \head{units} \\
\hline
 				 & \!\!\!\!\!\begin{tabular}{c|c} \head{sandstone}& \head{shale}\phantom{aaa} \end{tabular} &  \\ 
\hline
density ( $\rho^{s}$ ) & \begin{tabular}{c|c} 2648 & 2608 \end{tabular}&  kg $\mbox{m}^{-3}$  \\
\hline
specific thermal capacity ( $c^{s}$ ) & \begin{tabular}{c|c} 741 & 795 \end{tabular} & J $\mbox{K}^{-1}$ $\mbox{kg}^{-1}$  \\
\hline
specific heat conductivity ( $\lambda_s$ ) & \begin{tabular}{c|c} 3.0 & 1.2 \end{tabular} & W $\mbox{K}^{-1}$ $\mbox{m}^{-1}$  \\
\hline
sediment porosity at deposition ( $\phi_0$ ) & \begin{tabular}{c|c} 0.5 & 0.8 \end{tabular} & - \\
\hline
residual porosity ( $\phi_f$ ) & \begin{tabular}{c|c} 0.14 & 0.08 \end{tabular} & - \\
\hline
$k_1$, porosity-permeability law  & \begin{tabular}{c|c} 14.9 & 6.75 \end{tabular} & $\log(\mbox{m}^{2})$  \\
\hline
$k_2$, porosity-permeability law  & \begin{tabular}{c|c} \,\,\,\,\,\,\,\,\,\,\,\,\,1.94\, & uncertain \end{tabular} & $\log(\mbox{m}^{2})$  \\
\hline
\multicolumn{3}{|c|}{\textbf{QUARTZ CHARACTERISTICS}} \\  
\hline
\head{parameter} & \head{value} & \head{units} \\
\hline
density ( $\rho_Q$ ) & 2650 &  kg $\mbox{m}^{-3}$  \\
\hline
molar mass ( $M_Q$ ) & $6.008 \times 10^{-2}$ &  kg $\mbox{mol}^{-1}$  \\
\hline
initial specific surface area ( $A_0$ ) & $10^{4}$ &  $\mbox{m}^{-1}$  \\
\hline
parameter of reaction model ( $a_q$ ) & $ 5\times 10^{-19}$ &  mol $\mbox{m}^{-2}$  $\mbox{s}^{-1}$  \\
\hline
parameter of reaction model ( $b_q$ ) & 0.022 & $\mbox{Celsius}^{-1}$  \\
\hline
activation temperature ( $T_C$ ) & 373.15 &  K  \\
\hline
\end{tabular}
\end{center}
\caption{List of model parameters set to a fixed value. $k_2$ for shales is one of the parameters considered as uncertain, see Table \ref{table:uncertain-parameters}.}
\label{table:fixed-parameters}
\end{table}


We test the procedure described in Section \ref{sec:met} considering a synthetic test case characterized by alternating depositional events of two geomaterials, i.e. sand and shale sediments. 
Altough most of the effective parameters are affected by considerable uncertainty, given the large space-time scales involved in the compaction processes, here \Lc{we aim at investigating the characterization associated with input which have a large impact on the location of the interfaces between different materials}{we focus only on those parameters which have the strongest influence on variability of interfaces position.} With this idea in mind, we select three uncertain input parameters: $\beta_{sd}$ and $\beta_{sh}$, appearing in \EQU \eqref{eq:mec-comp}, i.e. the porous medium vertical compressibility for sandstone and shale, respectively; and the coefficient $k_2^{sh}$, which appears in \EQU \eqref{eq:k-phi} and defines the vertical permeability of shale. The ranges of variability assigned to these materials are listed in \TAB \ref{table:uncertain-parameters}, while \TAB \ref{table:fixed-parameters} reports the total sedimentation time, the sedimentation rate and the depositional order of sediments, together with boundary conditions and all other model parameters which are set to fixed values for this analysis.

We consider two cases characterized by two distinct parameter spaces: for case (A) only $\beta_{sd}$ and $\beta_{sh}$ are free to vary in selected variability ranges while in case (B) all three parameters are considered uncertain. We consider these two cases to quantify the effect of pure mechanical compaction on interfaces positions against the combined effect of mechanical compaction and vertical fluid flow inside the rock domain. 
In what follows, we build our sparse grids employing Gauss--Legendre nodes and a linear ``level-to-nodes'' function, $m(j)=j$,
with index sets \eqref{eq:iso-set} or \eqref{eq:aniso-set}.
\La{As already mentioned, all results presented hereafter refer to present time.}

\subsubsection{Approximation of interface positions}
\FIG \ref{image:compaction-histories} displays the distribution of geomaterials (shale and sandstone) obtained upon solving the full compaction model for the collocation points of a sparse grid sampling of the parameter space for the two cases. \Lc{Results refer to present time (i.e., the final configuration obtained for each simulation)}{} Each column in the two panels in the \FIG represents an independent simulation of the process: the juxtaposition of the results obtained for all collocation points allow to visually appreciate the variability of the interface position as a function of the uncertain parameters (left figure: isotropic sparse grid \EQU \eqref{eq:iso-set}, $w=6$, 137 collocation points; right figure: anisotropic sparse grid \EQU \eqref{eq:aniso-set}, $\aalpha_4=[4\,4\,1]$, $w=12$, 133 collocation points). 
The results show that interface $\Psi_1$ is constant in all realizations for both cases A and B. This is explained upon observing that the position of this interface is fixed and corresponds to the position of the basin top.  Interface $\Psi_2$ also display modest changes as a function of the model parameters for the two considered cases. The remaining interface locations $\Psi_3$-$\Psi_6$ display variations larger than 100m across the full sample of parameter realizations. We observe however  that for case A, for which $k_2^{sh}$ is fixed to a constant value, the location of interfaces $\Psi_3$-$\Psi_6$ change by maximum 300m. On the other hand when $k_2^{sh}$ is considered as uncertain (case B) the location of the deepest interfaces $\Psi_5$-$\Psi_6$ may vary by up to 1000m across the considered sample. 

\begin{figure}[tbp]
\begin{center}
\includegraphics[width=0.49\textwidth]{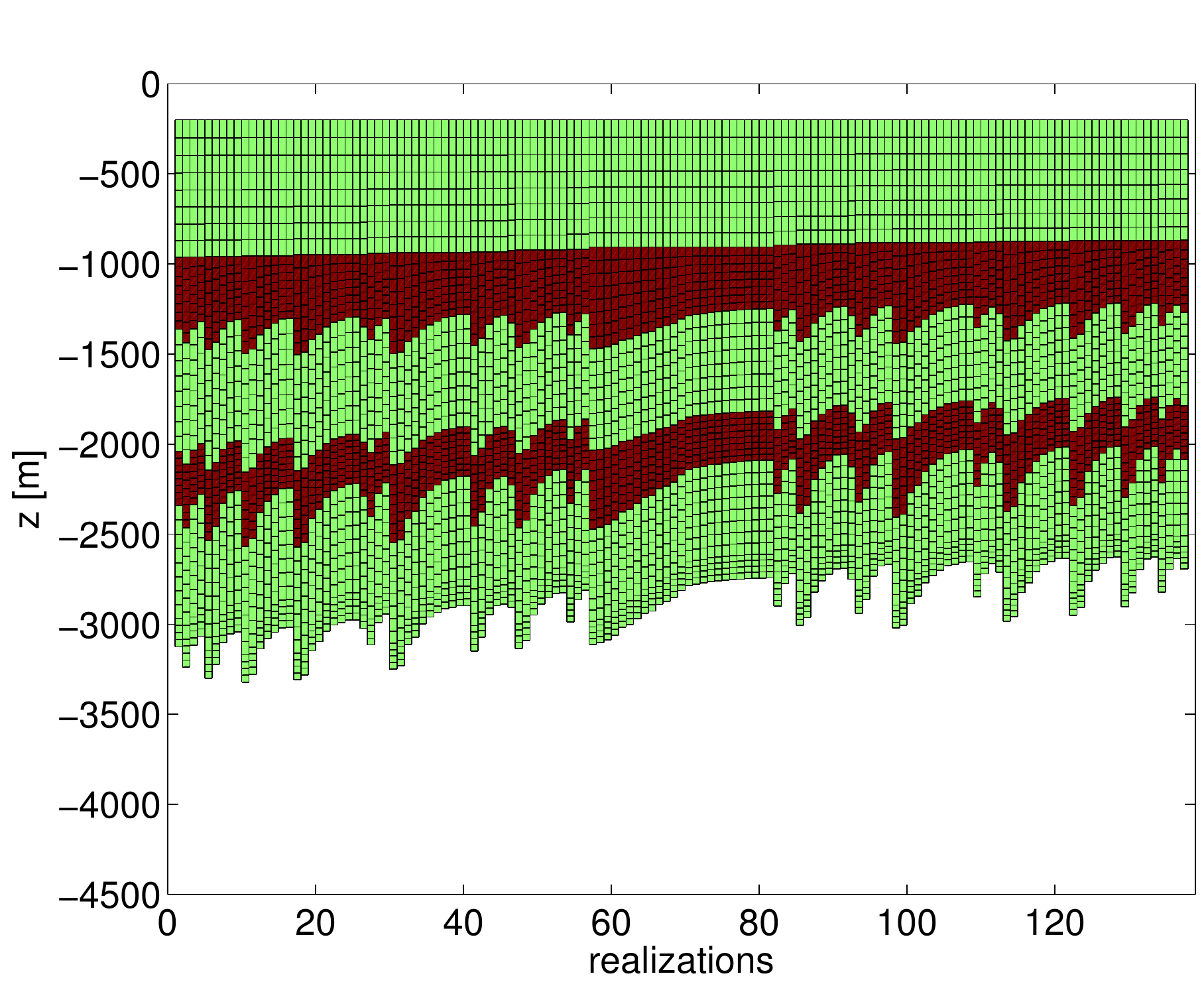}  
\includegraphics[width=0.49\textwidth]{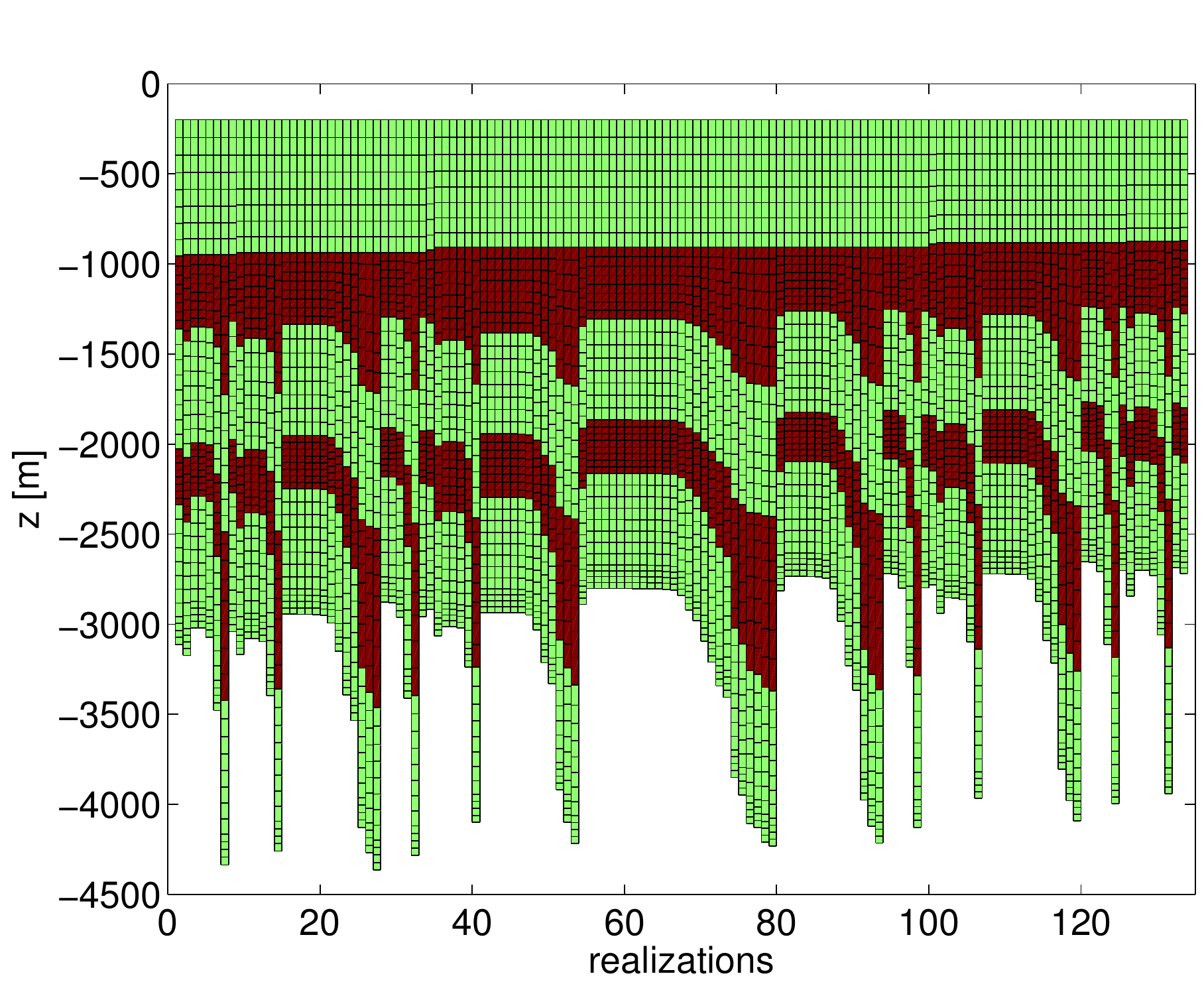}  
\caption{Compaction histories: list of full model realizations obtained for each set of input parameters. Results are referred to the analyzed cases: on the left (case A) only $\beta_{sd}$ and $\beta_{sh}$ are considered as uncertain, while on the right (case B) also $k_2^{sh}$  can vary in selected space interval.}
\label{image:compaction-histories}
\end{center}
\end{figure}

To identify the impact of each parameter on the variability of the interface position we resort to the total Sobol sensitivity indices $S^T$, which quantify the contribution of each uncertain parameter to the variance of the interface position (see \cite{feal:compgeo,lever.eal:inversion} for an efficient procedure to compute such indices starting from a sparse \Lc{grids}{grid} approximation). The results are depicted in \FIG \ref{image:Sobol} for case (A) and (B) respectively. Note that Sobol indices are not computed for the first interface since its position is known and fixed by the boundary condition.  

\FIG \ref{image:Sobol}-left shows the results obtained for case (A). 
Its analysis highlights that the position of each interface $\Psi_i$ is sensitive to compaction taking place within materials found at shallower locations than the interface itself. For example, we observe that the position of interface $\Psi_2$ is influenced merely by $\beta_{sd}$, sandstone being the only material overlying this interface. The parameter $\beta_{sh}$ influences the position of  $\Psi_3$ to a larger extent than $\beta_{sd}$, as this interface is directly underlying a shale lithological unit. For interfaces number $\Psi_4$, $\Psi_5$, $\Psi_6$ the influence of the two parameters tends to level off and the difference between $S^T_{\beta_{sd}}$ and $S^T_{\beta_{sh}}$ is reduced. 
These results suggest that for case (A) $\beta_{sd}$ and $\beta_{sh}$ have a similar impact on the uncertainty of the position of the interfaces.

\FIG \ref{image:Sobol}-right shows the results obtained for case (B). The role of $\beta_{sd}$ and $\beta_{sh}$ on the interface position display a similar behavior as in case (A). However, the influence of parameter $k_2^{sh}$ is predominant on the position of all interfaces with the exception of $\Psi_2$, which lies beneath the shallowest sandstone unit and therefore does not depend on the parameters governing flow and compaction in the shale layers. For all other interfaces $\Psi_i$ with $i = 3 \ldots 6$ we obtain $S^T_{k_2^{sh}} > 0.75$. This result indicates that the approximation of the material interfaces may benefit from anisotropic sampling with an increased sampling frequency along the parameter $k_2^{sh}$, upon following the criterion in \eqref{eq:aniso-set}.

\begin{figure}[tb]
\begin{center}
\includegraphics[width=0.5\textwidth,angle=0]{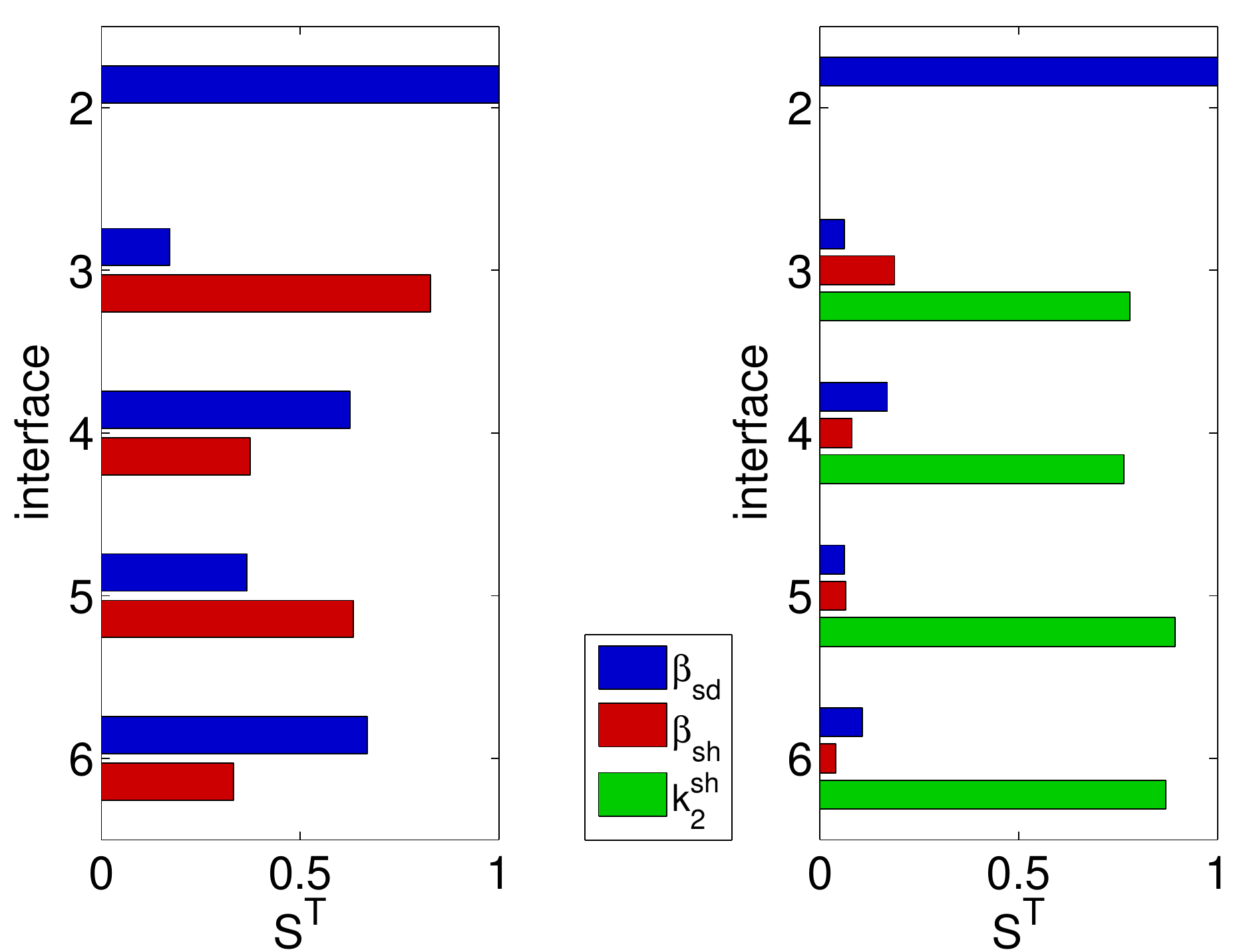} 
\caption{Total Sobol indices for uncertain parameters referred to estimated interface positions. Results are referred to the analyzed cases: on the left (A) only $\beta_{sd}$ and $\beta_{sh}$ are considered as uncertain and characterized by blue and red bars respectively, while on the right (B) also $k_2^{sh}$, green bars, can vary in selected space interval.}
\label{image:Sobol}
\end{center}
\end{figure}

In the following we refer to case (B) only, which displays larger variations in the interface positions and therefore is more challenging for our uncertainty quantification procedure.


We now aim at assessing the quality of the estimation of the interface positions $\Psi_k$ through sparse grids. 
To this end, we perform a convergence study of the sparse \Lc{grids}{grid} approximation of $\Psi_k$, measuring the error with
the two different error metrics detailed below. Moreover, we will repeat this study for
different types of sparse grids, 
by considering index sets as in \eqref{eq:aniso-set} with different choices of the vector $\aalpha$;
upon fixing the choice of $\aalpha$, the convergence study will be done increasing the index sets, 
hence the number of sparse grid points, by increasing $w$ in \eqref{eq:aniso-set}. 
The relation between the order $w$ and the number of sparse grid points, which we denote from here on as $N_{coll}$,
is depicted in \FIG \ref{image:w_Ncoll} for a selected set of sparse grids. 

\begin{figure}[tb]
\begin{center}
\includegraphics[width=0.4\textwidth,angle=0]{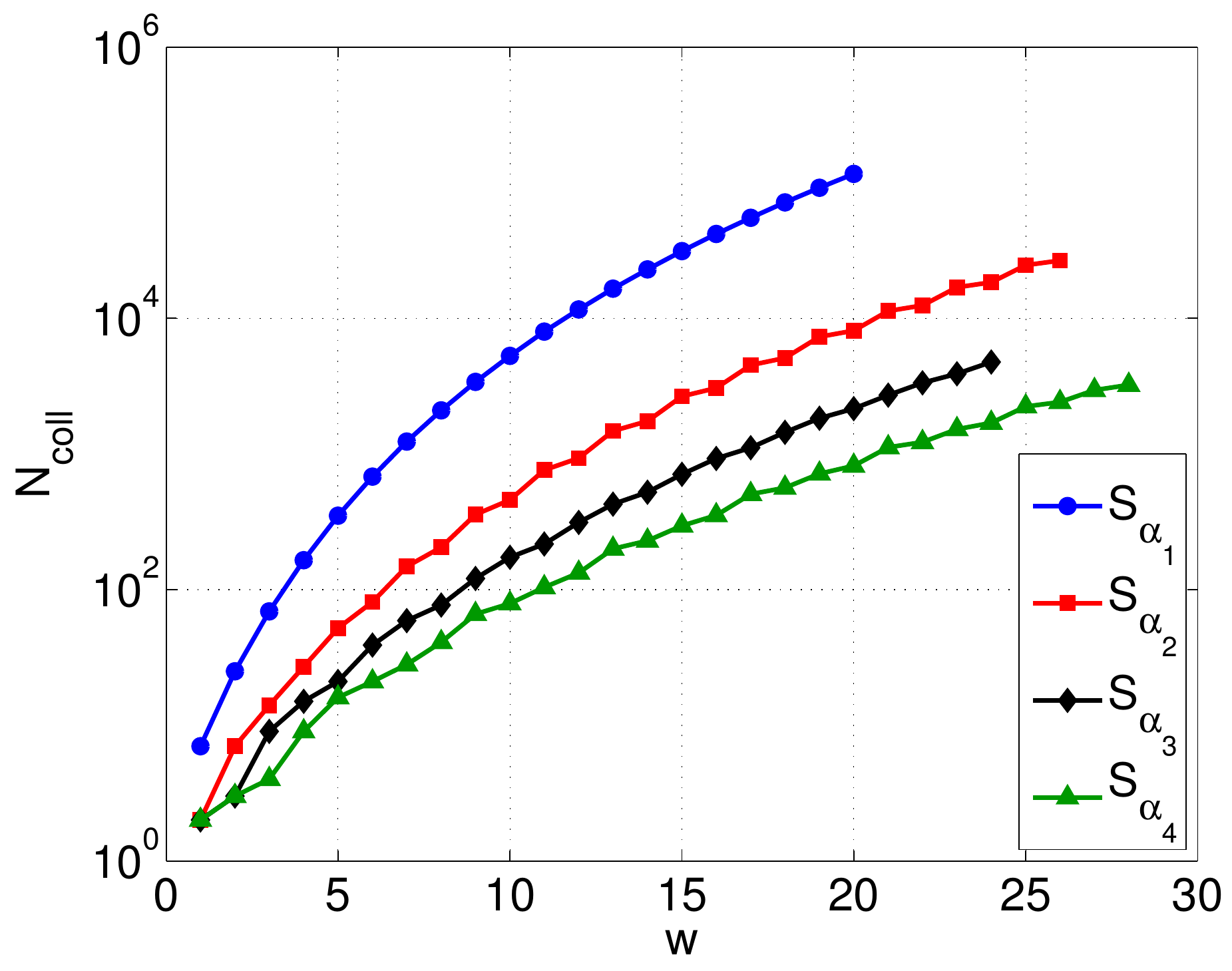} 
\caption{The relation between the order $w$ and $N_{coll}$ for a selected set of sparse grids:  $\aalpha_1=[1\,1\,1]$, $\aalpha_2=[2\,2\,1]$, $\aalpha_3=[3\,3\,1]$,
$\aalpha_4=[4\,4\,1]$.}
\label{image:w_Ncoll}
\end{center}
\end{figure}


Concerning the choice of $\aalpha$, we start from the observation 
anticipated by the results in \FIG \ref{image:Sobol}-right that the position of the interfaces
$\Psi_k$ is predominantly influenced by the parameter $k_2^{sh}$. 
Our choice is then to choose $\aalpha$ such that the resulting sparse grid
will be more refined along $k_2^{sh}$, while keeping the approximation degree identical
along $\beta_{sd}$ and $\beta_{sh}$. To do so, we consider the
following four choices for $\aalpha$: $\aalpha_1=[1\,1\,1]$, $\aalpha_2=[2\,2\,1]$, $\aalpha_3=[3\,3\,1]$,
$\aalpha_4=[4\,4\,1]$. Clearly, choosing $\aalpha=\aalpha_1$ in \eqref{eq:aniso-set} is equivalent to building
an isotropic sparse grid, i.e. using the set defined in \eqref{eq:iso-set}, in which all parameters are identically refined. 
According to the notation introduced above, the grids thus obtained should be 
denoted by $\mcS_{\mcI_{\aalpha_i}}$, with $i=1,\ldots,4$; however, we will employ a lighter notation, i.e.
$\mcS_{\aalpha_i}$. Similarly, the associated quadrature rules will be denoted by $\mcQ_{\aalpha_i}$, cf. equation \eqref{eq:sparse-quad}. 

We employ the following two metrics to evaluate the accuracy of the various sparse grid approximations:
\begin{itemize}
\item the relative error of the mean position obtained through the sparse approximation 
\begin{equation} \label{eq:mean-error}
E_{mean}\left(\Psi_k\right) = \frac{\left| \mcQ_{\aalpha_i}\left(\Psi_k\right) - \mu_{ref}\left(\Psi_k\right)\right|}{\left| \mu_{ref} \left(\Psi_k\right) \right|}
\end{equation}
where $k = 1, \ldots ,6$ identifies the interfaces and $\mu_{ref} \left(\Psi_k\right)$ is a reference value for the mean of the interface position.
We consider here as reference value the mean obtained with a very refined version of the sparse grid 
$\mcS_{\mcI_{\aalpha_4}}$, with around 30000 collocation  points.

\item the maximum norm of the error of the interface position obtained with respect to the one predicted by the full model simulation within the parameter space 
\begin{equation} \label{eq:norma-error}
E_{max}\left(\Psi_k \right) = \left\lVert \frac{\left[\mathcal{S}_{\aalpha_i}\left[\Psi_k\right] \left(\mathbf{p}\right) - \Psi_k^{FM}\left(\mathbf{p}\right)\right]}{\Psi_k^{FM}\left(\mathbf{p}\right)} \right\rVert _\infty
\end{equation}
where $\Psi_k^{FM}\left(\mathbf{p}\right)$ is the value of $\Psi_k$ obtained from the direct numerical approximation of the full model. The metric \eqref{eq:norma-error} is numerically computed by employing a random sampling of the functions $\mathcal{S}_{\aalpha_i}\left[\Psi_k\right] \left(\mathbf{p}\right)$ and $\Psi_k^{FM}\left(\mathbf{p}\right)$  for $\mathbf{p } \in \Gamma$. We consider here 1000 random evaluations to compute \eqref{eq:norma-error}. 
\end{itemize}

\FIG \ref{image:convergences} displays the two metrics $E_{max}$ and $E_{mean}$ as a function of the number of realizations required to build each sparse grid approximation. We consider here the variation of the two metrics for the positions $\Psi_2$ and $\Psi_5$. Note that Figure \ref{image:convergences} also reports the results on the convergence of the mean value of $\Psi_2$ and $\Psi_5$ computed by a straightforward Monte Carlo approximation. For interface positions $\Psi_3, \Psi_4, \Psi_6$  the results are qualitatively similar to those obtained for $\Psi_5$ and for this reason are not reported here. 

Results reporting the convergence of interface $\Psi_2$ are reported in \FIG \ref{image:convergences}a-b. A first observation emerges from the comparison between sparse grid and Monte Carlo sampling technique, where the more efficient behavior of the former is highlighted. For example, we find $E_{mean} \approx 10^{-8}$  for a number of collocaction points within 60 and 200 depending on the type of sparse grid employed. For the same number of realizations, the mean error associated with the Monte Carlo approach is close to $10^{-3}$. The position of $\Psi_2$ is approximated with high accuracy by sparse grid surrogate models, with $E_{mean} \approx 10^{-14}$ for $N_{coll} \approx 10^3$ (see \FIG \ref{image:convergences}a) and $E_{max} \approx 10^{-12}$ (see \FIG \ref{image:convergences}b). Among the sparse grid approximations $\mathcal{S}_{\aalpha_3}$, $\mathcal{S}_{\aalpha_4}$ are the least effective ones. This can be explained upon observing that $\Psi_2$ only depends on $\beta_{sd}$  (\FIG \ref{image:Sobol}b), so increasing the sampling frequency along $k_2^{sh}$ adds to the sparse grid points that are not improving
the quality of the approximation.

\FIG \ref{image:convergences}c-d displays the results obtained for $\Psi_5$. 
 For this interface, 
 \Lc{the sparse grid approximations outperform the Monte Carlo approximation 
   by a factor comprised between 10 and 100 for a fixed number of sampling points,}{
   the anisotropic sparse grids approximations
   $\mathcal{S}_{\aalpha_3}$ and $\mathcal{S}_{\aalpha_4}$ still
   significantly outperform the Monte Carlo approximation, both in
   terms of error size and convergence rate} 
 when metric $E_{mean}$ is considered (\FIG \ref{image:convergences}c) 
 \La{although the improvement is in this case less than for the second interface, 
   due to the increased complexity of the function to approximate, 
   see Figures \ref{fig:smooth-and-non-smooth} and \ref{image:compaction-histories}}.
\Lc{Moreover, the anisotropic sparse grid approximation
  $\mathcal{S}_{\aalpha_4}$ converges faster than the others (see \FIG
  \ref{image:convergences}c-d), while in particular the isotropic
  sparse grid $\mathcal{S}_{\aalpha_1}$ is the least efficient choice
  (as expected).}{}
This result shows that anisotropic sampling of the parameter space
yields a computational advantage if compared on isotropic one for the
approximation of deeply buried interfaces.  In the remainder of this
section, to test the applicability of the sparse grid technique we
select the sparse grid $\mathcal{S}_{\aalpha_4}$ and fix $ w =12$,
which yields $N_{coll} =133$.

\begin{figure}[tb]
  \begin{subfigure}{0.5\textwidth}
    \centering
    \includegraphics[width=0.8\linewidth]{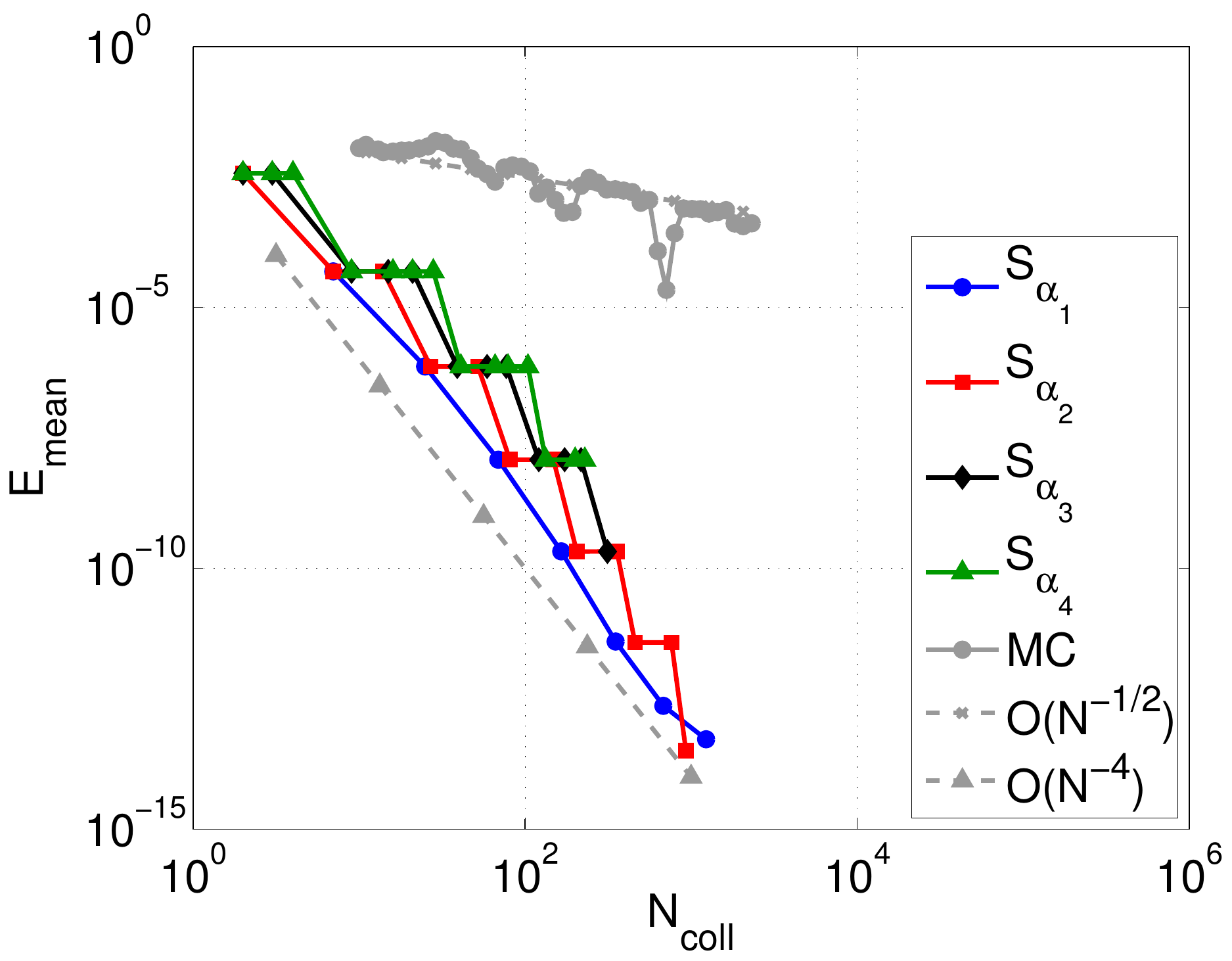} 
    \caption{mean trends for second interface}
    \label{image:convergence-mean-2}
  \end{subfigure}
  \begin{subfigure}{0.5\textwidth}
    \centering
    \includegraphics[width=0.8\linewidth]{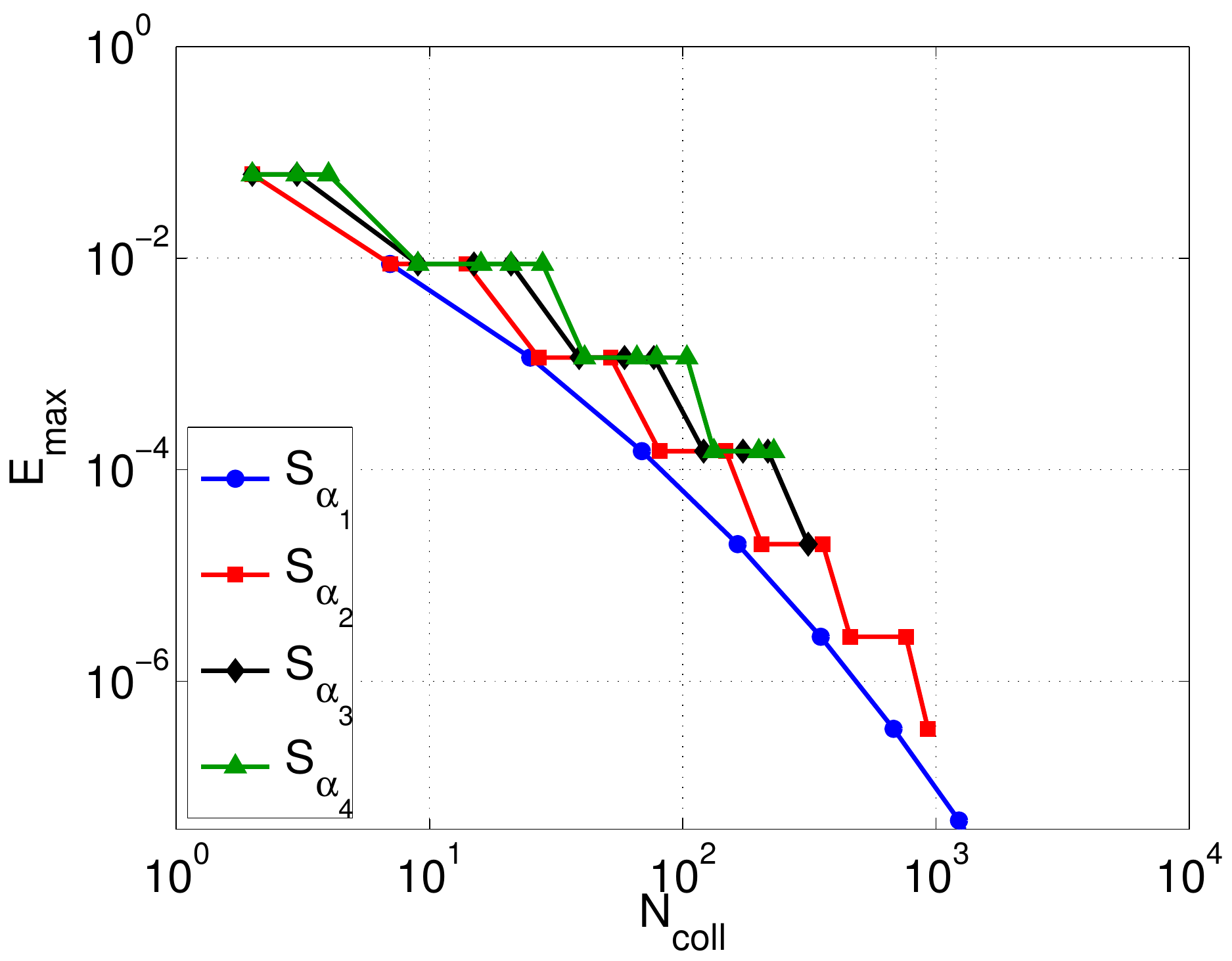}
    \caption{maximum norm trends for second interface}
    \label{image:convergence-norm-2}
  \end{subfigure}
  \begin{subfigure}{0.5\textwidth}
    \centering
    \includegraphics[width=0.8\linewidth]{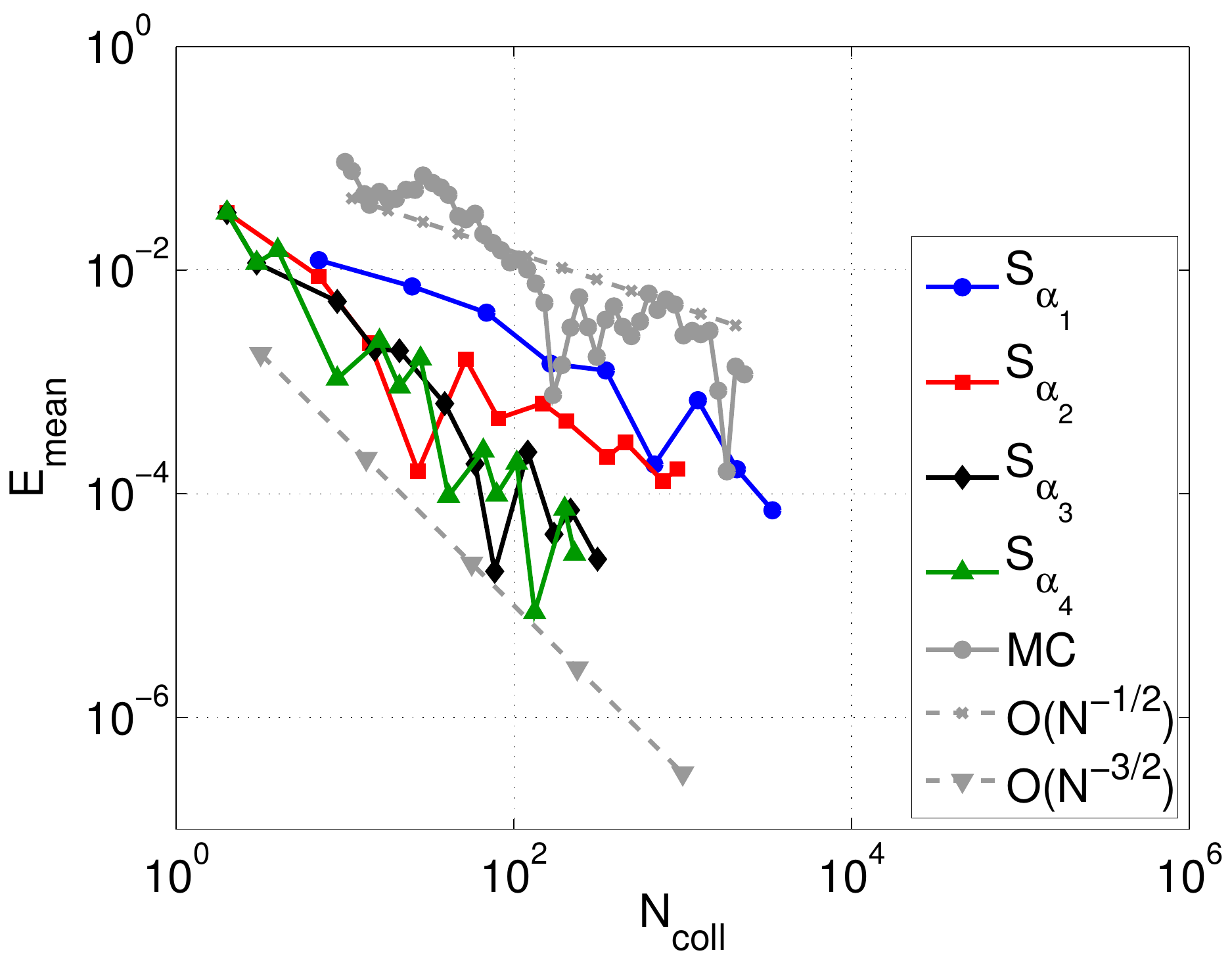} 
    \caption{mean trends for fifth interface}
    \label{image:convergence-mean-5}
  \end{subfigure}
  \begin{subfigure}{0.5\textwidth}
    \centering
    \includegraphics[width=0.8\linewidth]{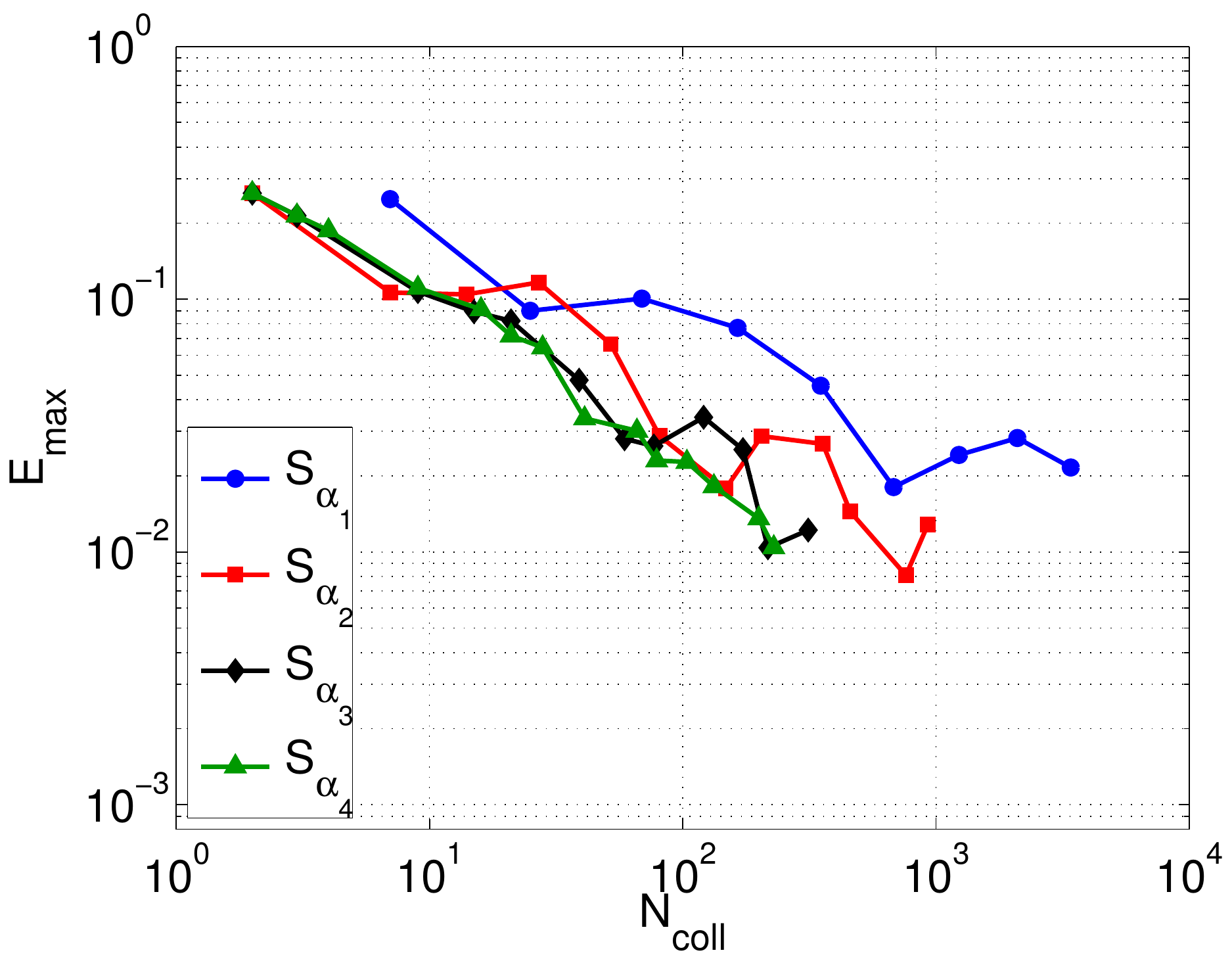}
    \caption{maximum norm for fifth interface}
    \label{image:convergence-norm-5}
  \end{subfigure}
  \caption{Study of the convergence of mean and maximum norm}
  \label{image:convergences} 
\end{figure}

\Lc{Before going further deep in the analysis of the results
we also briefly comment on the computational gain offered by the sparse grids approximation.
On a processor with clock-speed 3.30GHz and 15 MB cache, a full model run on Matlab 8.5
takes about 2 seconds, depending on the value of the random parameters.
Therefore, comparing the results of Monte Carlo and sparse grids approximations in Figure 
\ref{image:convergence-mean-5}, an approximation on the mean position
of the fifth interface with a $0.1\%$ accuracy takes about $2000 \times 2 = 4000$ seconds 
(roughly 1 hour) while sparse grids need at most
$300 \times 2 = 600$ seconds (10 minutes) with an isotropic sparse grid, 
and around $50\times 2=100$ seconds (less than 2 minutes) for the anisotropic sparse
grids (the time needed to actually determining the set of evaluation points of the sparse grid is negligible). 
Moreover, upon obtaining the samples needed to build the sparse grid, 
evaluating the position of an interface for a new value of the random parameters is extremely efficient,  
in the order of a fraction of a millisecond per evaluation; 
the same holds for evaluating the porosity at a given depth with the two-steps procedure
described in Section \ref{sec:disc-output}. 
Thus, the analyses that will be performed below, that need to evaluate the sparse
grid approximation over a Monte Carlo ensamble with a few thousand samples, can be performed in a few seconds.
We also remark that the sparse grid evaluation algorithm that we use, see \cite{lorenzo:sparse-grids-code},
does not scale linearly with the number of points to be evaluated nor with the number of points of the 
sparse grids, due to some internal optimizations,
therefore it would be meaningless to provide a precise value of the time needed by a single evaluation.}{}
 

The estimation of the material interface positions $\Psi_i$ allows to
estimate the probability to find a specific material at a
selected depth in the system. We deal with a conceptual model
characterized by alternating deposition of different materials which
divide the domain into five macro layers, as listed in the top part of
\TAB \ref{table:uncertain-parameters}. The probabilistic distribution
of geomaterials along depth can be highly valuable when performing
history matching of real sedimentary systems. To obtain this result,
we consider now a series of $N_z$ points between $\tilde{z}_k \in
[-500,-2300]$m separated by a constant interval of 10m. We then
evaluate which material $M(\tilde{z}_k)$ occurs at each of these
locations, upon comparing $\tilde{z}_k$ with the estimated interface
locations $\Psi_i$. Note that since $\Psi_i$ depends on the random
parameters ${\bf p}$, then $M(\tilde{z}_k)$ is also random.  \FIG
\ref{image:material-occurrence} shows the vertical distribution
of probability obtained upon computing the relative frequency
$Freq(M(\tilde{z}_k) = M_i)$, i.e. the probability of observing
$M(\tilde{z}_k) = M_i$ with $M_i = 1, \ldots ,5$. To compute
$Freq(M(\tilde{z}_k) = M_i)$ we evaluate the sparse grid approximation
over a set of 5000 random points. 
Figure \ref{image:material-occurrence}-left shows the vertical distribution of
$Freq(M(\tilde{z}_k) = M_i)$ evaluated through the selected sparse
grid approximation. Note that for $z > -1700$m a single depth is
associated with nonzero probability to observe a single material or
two materials, close to the interface. For depths deeper than 1700
meters the complexity increases and a single location may be
associated with up to three different materials, i.e. $M_3$, $M_4$ and
$M_5$. These results provide a quantitative assessment of the
qualitative observation emerging from \FIG \ref{image:compaction-histories}b.

In \FIG \ref{image:material-occurrence}-right we complement the results in \FIG \ref{image:material-occurrence}-left upon showing the relative frequency of misclassified realizations ($Freq_{mis}$), i.e. of realizations for which the value of $M(\tilde{z}_k)$ estimated through the surrogate sparse grid model is different from the one obtained through the direct numerical simulation of the full model for the same parameter combination. The number of misclassified points increases with depth and with the probability to find more than one material at the same depth: in those regions in which the extent of overlapping is small, i.e. between $-850<z<950$m (i.e. for locations close to $\Psi_2$), the frequency is lower than 0.1\% while for the deepest locations considered, i.e. for $2000<z<2300$m, the relative frequency of misclassified point increases up a maximum value of 0.8\%. 
This value confirms that the implemented methodology provides a robust estimation of the material interface positions, and can be effectively employed in practice to obtain the probabilistic estimation of the vertical extent associated with all geological units. 

\begin{figure}[tbp]
\begin{center}
\includegraphics[width=0.55\textwidth,angle=0]{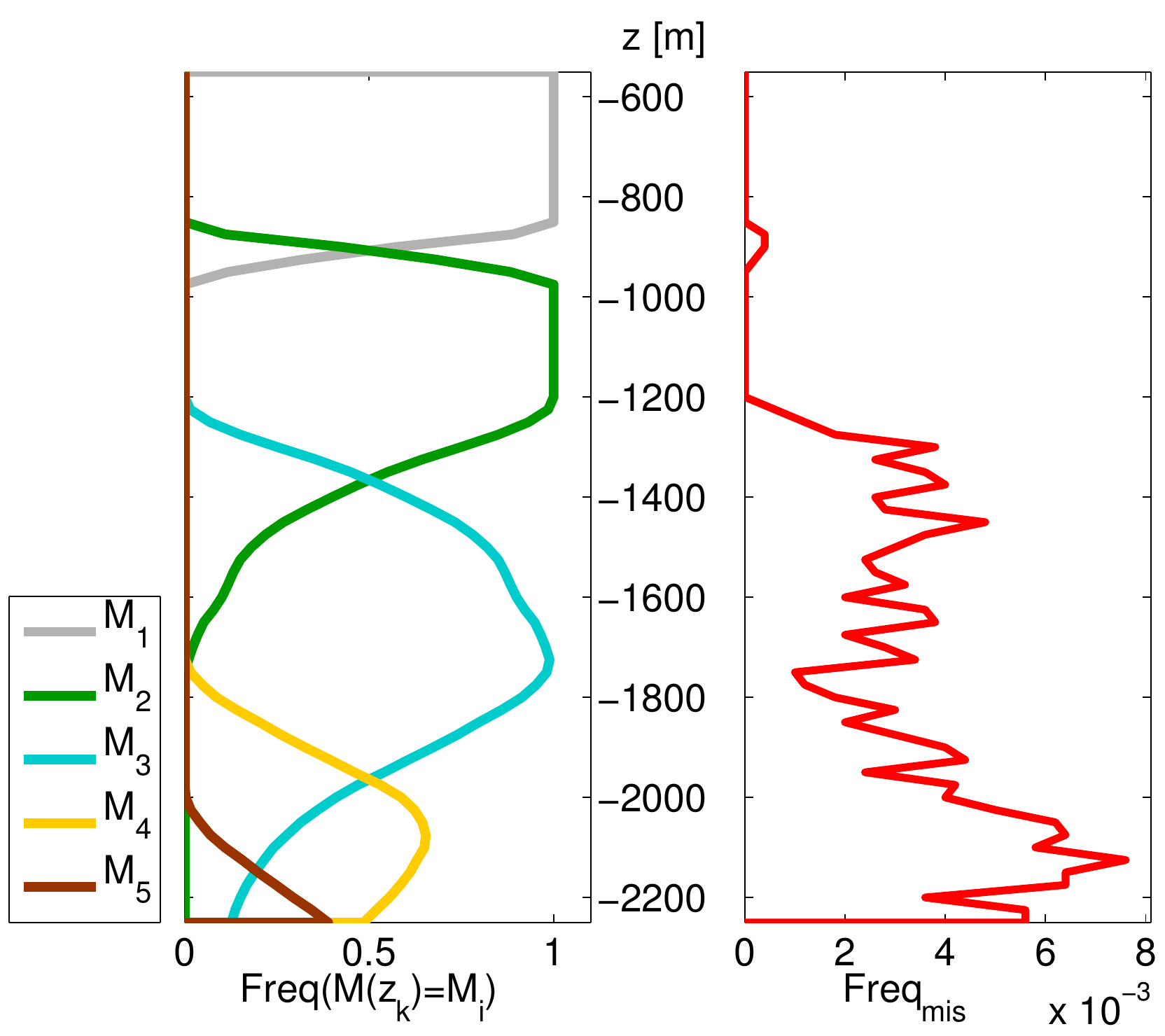}  
\caption{Left: frequency of material occurrences along the vertical direction. Right: frequency of misclassified materials 
  along the vertical direction}
\label{image:material-occurrence}
\end{center}
\end{figure}

\subsubsection{Uncertainty quantification for porosity distributions}

Finally, we investigate the propagation of uncertainty from input parameters to model outputs considering the probability density functions of porosity at different depths. In order to validate our procedure, the results of the two-step surrogate model procedure are compared against full model results in a Monte Carlo framework. 
Specifically, we consider a subset of the locations $\tilde{z}_k$ introduced above, we compute for each of them
an approximation of the \La{porosity} pdf by sampling our sparse grid surrogate model, and we validate the results upon comparing them with those rendered by the full model sampling. As remarked earlier, porosity as a state variable is heavily influenced by the characteristic of the geomaterials.

\La{To assess the accuracy of the results given by the two-step procedure we compare the results in terms of the distance between the cumulative distribution functions (cdfs) of porosity obtained across a range of depths $\tilde{z}_k \in [-600 -2200]$ m. This is evaluated as 
\begin{equation}
D = \max\left[ \left|F_{FM}\left(\phi\right) - F_{SG}\left(\phi\right)\right|\right]
\label{eq:distance}
\end{equation}
where $F_{FM}(\phi)$ and $F_{SG}(\phi)$ represent the cdfs of porosity
obtained through the numerical solution of the full model or by the
two-step model reduction technique introduced in Section
\ref{sec:disc-output}, respectively.  The left panel of Figure
\ref{image:ks-distance} displays the variation of $D$ along the
vertical direction. The results displayed here are obtained upon
approximating $F_{FM}(\phi)$ and $F_{SG}(\phi)$ by 5000 Monte Carlo
realizations. We consider here various candidate sparse grids
surrogate models, which are all obtained from sparse grid
approximations $\mathcal{S}_{\aalpha_4}$ upon fixing $w = 2, 6, 12$. 
We observe that $D$ is reduced to values of the order of
$10^{-2}$ for $w = 12$, 
while values typically larger
than $10^{-1}$ are attained for $w = 2$. The center and right panels
of Figure \ref{image:ks-distance} display the cdfs of porosity for two
depths, i.e. at $\tilde{z}_k = -900$ m and $\tilde{z}_k = -2000$ m. 
We select in particular the depth $\tilde{z}_k = -900$ m as it is
associated with the largest errors for $w = 12$. 
We observe that at both investigated depths the sparse grid surrogate 
associated with $w = 12$ 
approximates the full model response up to a reasonable accuracy, 
while large errors are observed for $w = 2$. 
This result is consistent with the observation that the metric \eqref{eq:distance} is influenced
by the accuracy achieved in approximating the location of interfaces, which is discussed above.  }

\begin{figure}[tbp]
\begin{center}
\includegraphics[width=1.0\textwidth,angle=0]{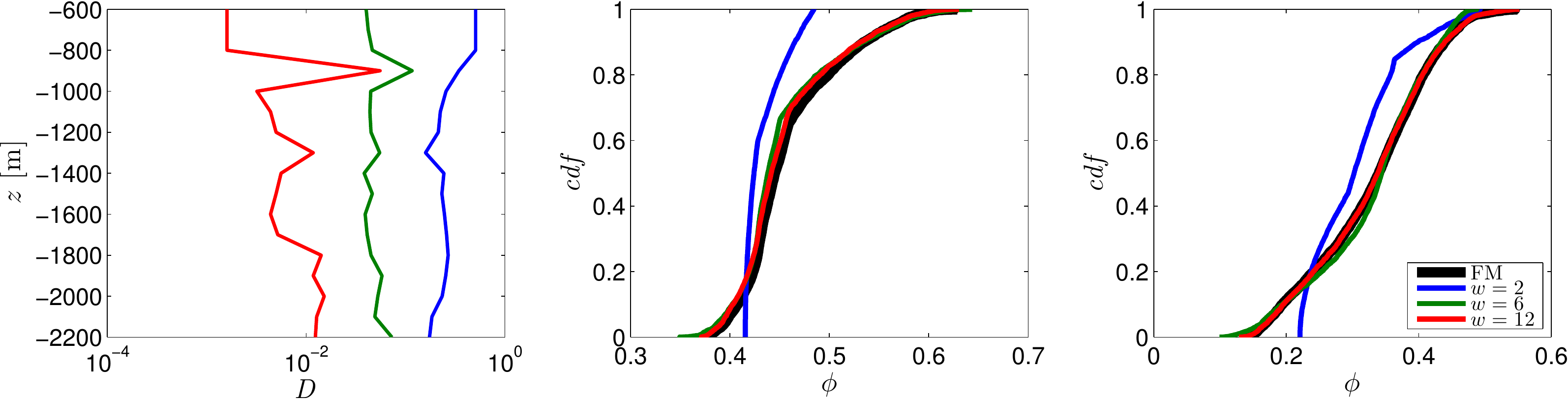}  
\caption{Numerical convergence of two-step surrogate model: variation of the distance $D$ \eqref{eq:distance} along the vertical direction 
for $\mathcal{S}_{\aalpha_4}$ with $w = 2, 6, 12$ (left), 
and comparison of sample cdfs of porosity obtained at depths  $\tilde{z}_k = -900$ m (center) and  $\tilde{z}_k = -2000$ m (right).} 
\label{image:ks-distance} 
\end{center}
\end{figure}

We showcase finally the capabilities of our technique upon analyzing in close details the full pdfs of porosity for some specific depth. Since the occurrence of up to three different materials can be observed at some locations, we expect to observe multimodal pdfs of porosity values for those  $\tilde{z}_k$. \La{The pdfs are obtained with Matlab command \texttt{ksdensity},
which is based on a normal kernel function; we have tuned manually the bandwidth of the kernel smoothing window to 0.02, 
to highlight the multimodality of the pdf.}

\FIG \ref{image:pdf-comparison} considers four depths to examine the shape of \La{porosity} pdfs obtained with both the sparse grid and the full model solution. In order to sample the vertical domain in a comprehensive way we show here the results associated with four locations, one from the upper part (-600 m), two from the middle (-1350 m, -1500 m) and one from the deepest part of the basin (-2250 m). In addition to the pdfs of porosity we also display for each selected location a scatterplot which allows for a direct comparison between the porosity values predicted with the sparse grid approximation and those given by the full model for the same location and for all considered Monte Carlo realizations of the parameters. In each scatterplot we distinguish the points associated with different materials $M_i$ and we also explicitly tag with red markers the misclassified points introduced above. As a preliminary remark, we observe that the porosity values found for some of the Monte Carlo realizations attain high values (up to 0.8) which are unlikely to be found in a real sedimentary setting (see Figure \ref{image:pdf-comparison}). Here, we detail the results of the full analysis which may be interpreted as a preliminary uncertainty analysis where large uncertainty bounds are assigned to model parameters. In principle, our model reduction methodology allows then to efficiently identify\Lc{ing }{} the regions of the parameter space where those values are found and focus on a restricted set of the full parameter space upon excluding unphysical or unrealistic combinations. These may be identified for example from an expected porosity trend which may come from direct field observation or other prior knowledge on the system.

As a first general observation the pdfs yielded by the sparse grid approximation are in close agreement with those obtained through the direct simulation of the full model. Some small differences are only observed for cases where the pdfs shape have a particularly complex and multimodal shape  (Figure \ref{image:pdf-comparison1350}-\ref{image:pdf-comparison2250} ).  
At $\tilde{z} = -600$m we find material $M_1$ (most recent sandstone layer), for all the Monte Carlo realizations consistently with results in \FIG \ref{image:material-occurrence}. \FIG \ref{image:pdf-comparison600} and \FIG \ref{image:scatterplot-600} show that the values of porosity obtained with the full model and the surrogate basically coincide and span the interval between 0.35 and 0.5. 
Location $\tilde{z} = -1350$m is associated with similar probability of observing material $M_2$ and $M_3$, see \FIG \ref{image:material-occurrence}a. The shape of the pdf in \FIG \ref{image:pdf-comparison1350} shows three main peaks, which can be tied to the porosity values attained within the two materials which are found at this location across the complete set of realizations of the Monte Carlo sample. This result is therefore informative on the compaction regime which may be observed at this location as a function of the assumed uncertainty of the selected random parameters. In particular we observe that: 
\begin{itemize}
\item The scatterplot in Figure \ref{image:scatterplot-1350} puts in evidence that realizations for which material $M_3$ (sandstone) is found at this location tend to cluster around $\phi \approx 0.3$, this value corresponding to the location of a peak of the porosity pdf (see Figure \ref{image:pdf-comparison1350}). 
\item Realizations for which material $M_2$ (shale) is found attain values of porosity of approximately 0.75 (see Figure \ref{image:scatterplot-1350}), which correspond to a location of a second peak of the porosity distribution (see Figure \ref{image:pdf-comparison1350}). These high values of porosity are explained upon observing that the pressure gradient in shale material may become larger than hydrostatic when shale permeability decreases (i.e. as a function of the value assumed by the random parameter $k_2^{sh}$). In these conditions the vertical effective stress decreases and this has a direct effect on porosity through equation \eqref{eq:mec-comp} (see e.g., \cite{Colombo16_HJ} for a complete discussion of this phenomenon).
\item a third peak is observed in the porosity pdf for $\phi \approx 0.6$, this value corresponding to the overlap of the porosity distributions in the two materials  $M_2$ and $M_3$. Note that a few misclassified points appear for $\phi \approx 0.5$, i.e. for porosity values which are possibly found in both materials.
\end{itemize}

At $\tilde{z} = -1500$m materials $M_2$ and $M_3$ may be found, with $Freq[M(\tilde{z}_k = -1500) = 3] > Freq[M(\tilde{z}_k = -1500) = 2]$ (Figure \ref{image:material-occurrence}a), i.e. larger probability associated with the presence of $M_3$ (sandstone) than $M_2$. \Lc{as a result the}{As a result, the } porosity pdf (\FIG \ref{image:pdf-comparison1500}) has here two distinct peaks which are found at porosity values consistent with those observed at depth $\tilde{z} = -1350$m, i.e. $\phi \approx 0.78$ corresponding to shale material and $\phi \approx 0.3$ for sandstone material. 

Finally, \FIG \ref{image:pdf-comparison2250} and \FIG \ref{image:scatterplot-2250} display the results obtained for location $\tilde{z} = -2250$m. At this location three different lithological units may be found as a function of the values assumed by the three random parameters across the Monte Carlo sample (see \FIG \ref{image:material-occurrence}a): the largest probability of occurrence is associated with material $M_4$ (shale) followed by $M_5$(sandstone) and $M_3$(sandstone). In this case the scatter plot shows a slight increase of the dispersion of the points around the bisector, which indicates that the accuracy of the sparse grid approximation is lower than that observed for shallower locations (compare Figure \ref{image:scatterplot-2250} e.g.  with  \ref{image:scatterplot-600}). The comparison of the pdfs displayed in Figure \ref{image:pdf-comparison2250} evidences however that these inaccuracies have only a minor consequence on the estimation of the pdf of $\phi$. Once again, we can tie the observed probability density of $\phi$ to the occurrence of the various geomaterials and relate this to the compaction behavior. The sandstone materials ($M_3$ and $M_5$) tend both to attain smaller porosity values than the shale unit. Note that as the burial depth increases, the maximum observed porosity decreases, i.e. we find here $ \phi < 0.7$ also within shale material $M_4$.

\begin{figure}[tbp] 
  \begin{subfigure}{0.5\textwidth}
    \centering
    \includegraphics[height=5cm]{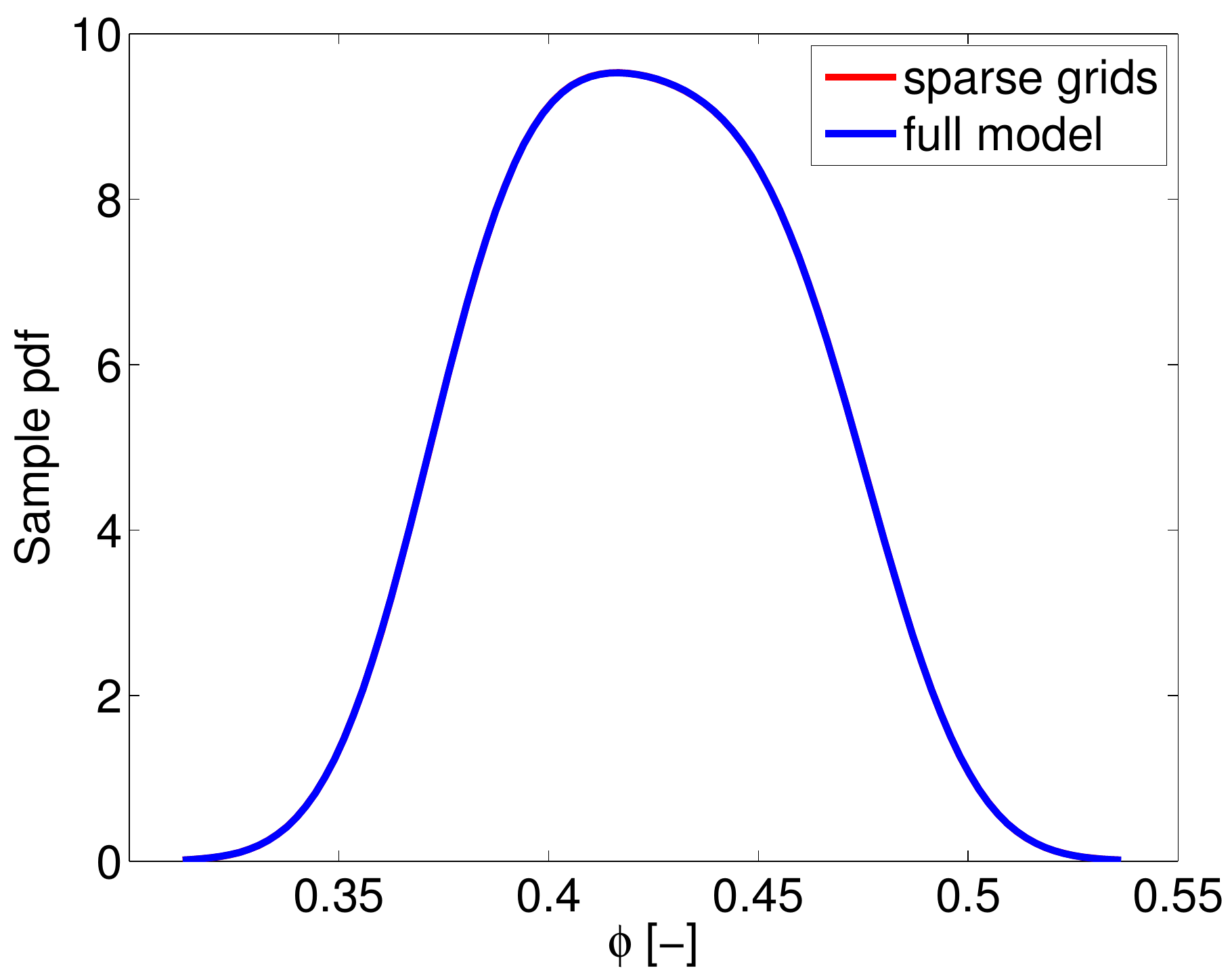} 
    \caption{pdf at -600 meters}
    \label{image:pdf-comparison600}
  \end{subfigure}
  \begin{subfigure}{0.5\textwidth}
    \centering
    \includegraphics[height=5cm]{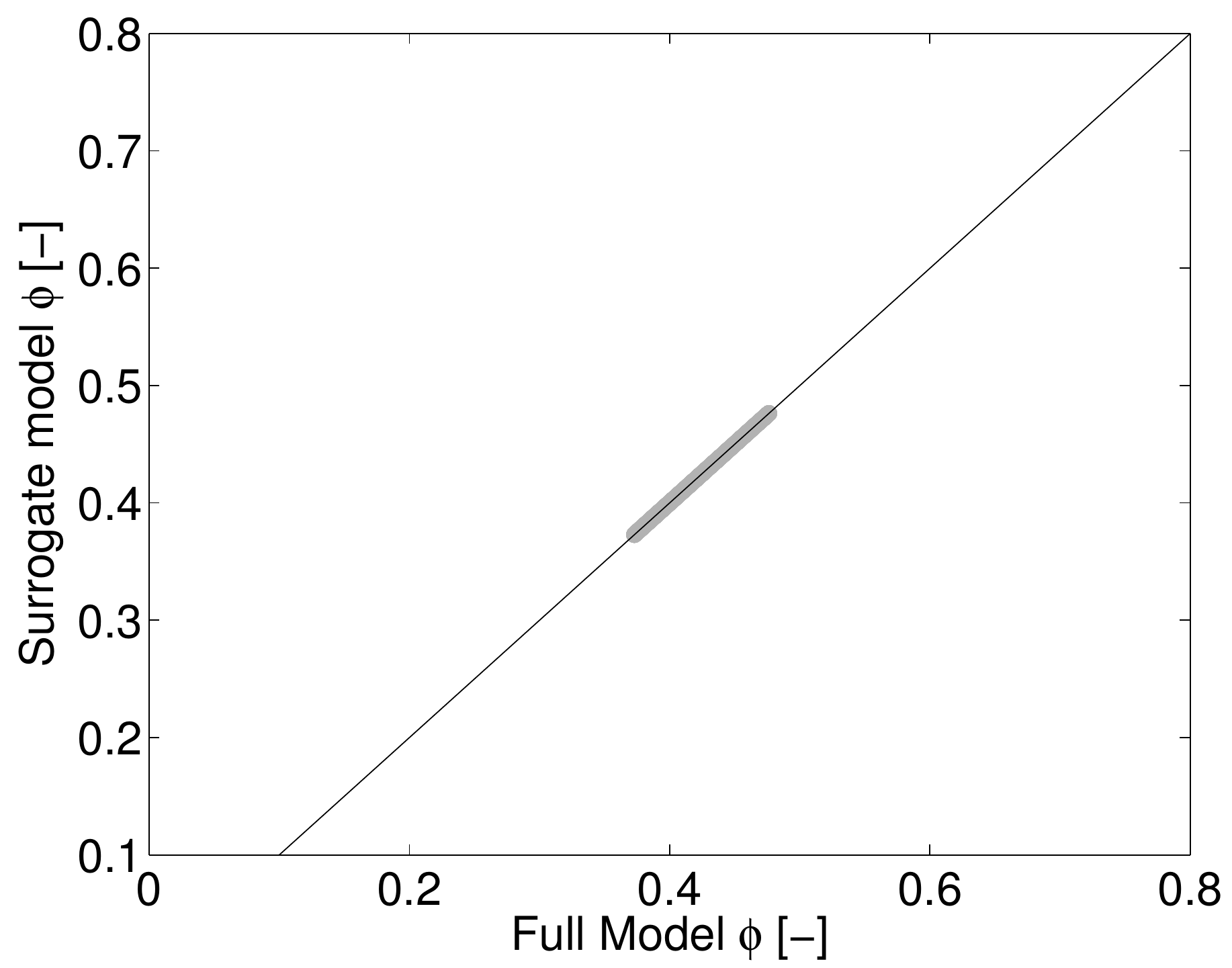}
    \caption{scatterplot at -600 meters}
    \label{image:scatterplot-600}
  \end{subfigure}
  \begin{subfigure}{0.5\textwidth}
    \centering
    \includegraphics[height=5cm]{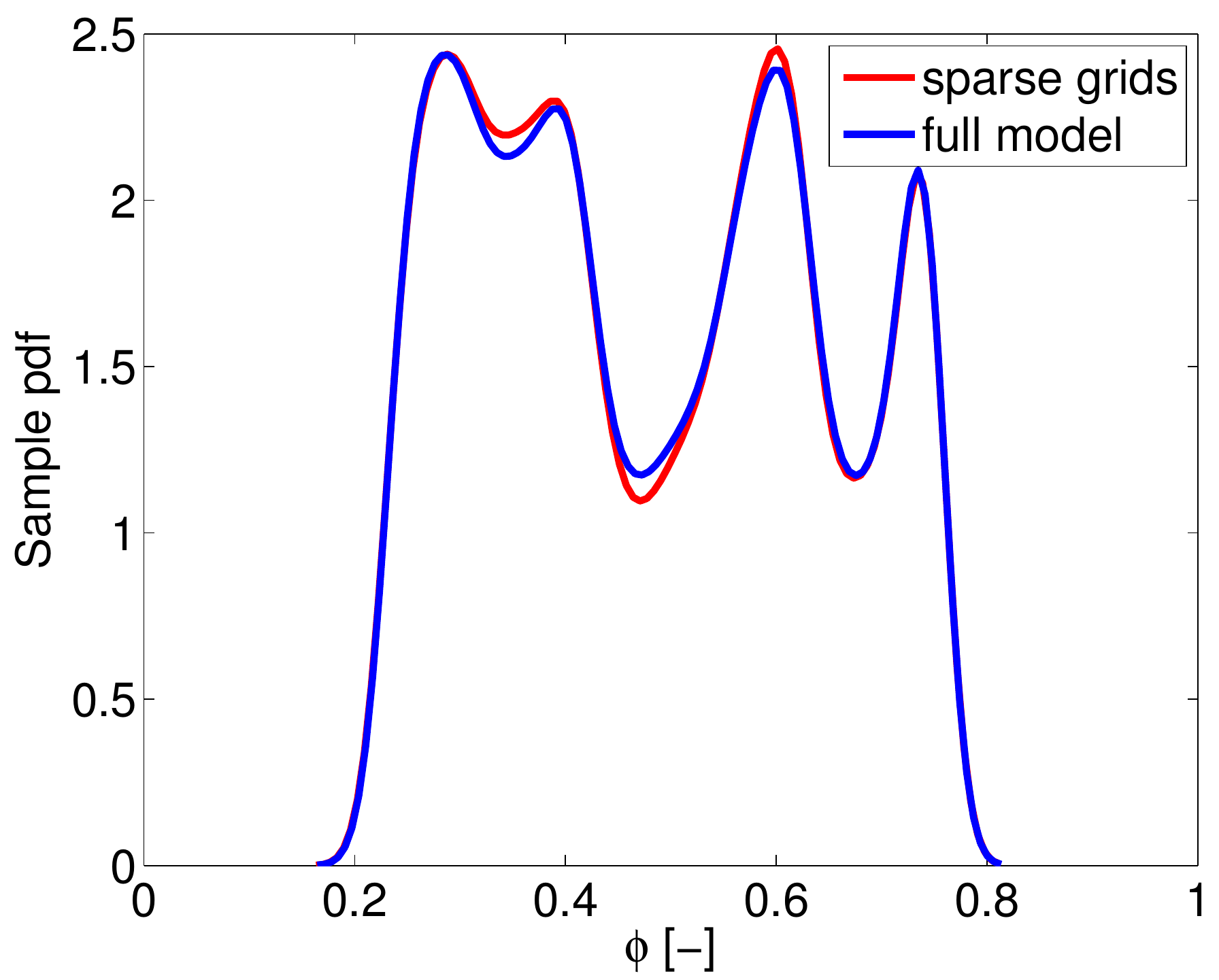}
    \caption{pdf at -1350 meters}
    \label{image:pdf-comparison1350}
  \end{subfigure}
  \begin{subfigure}{0.5\textwidth}
    \centering
    \includegraphics[height=5cm]{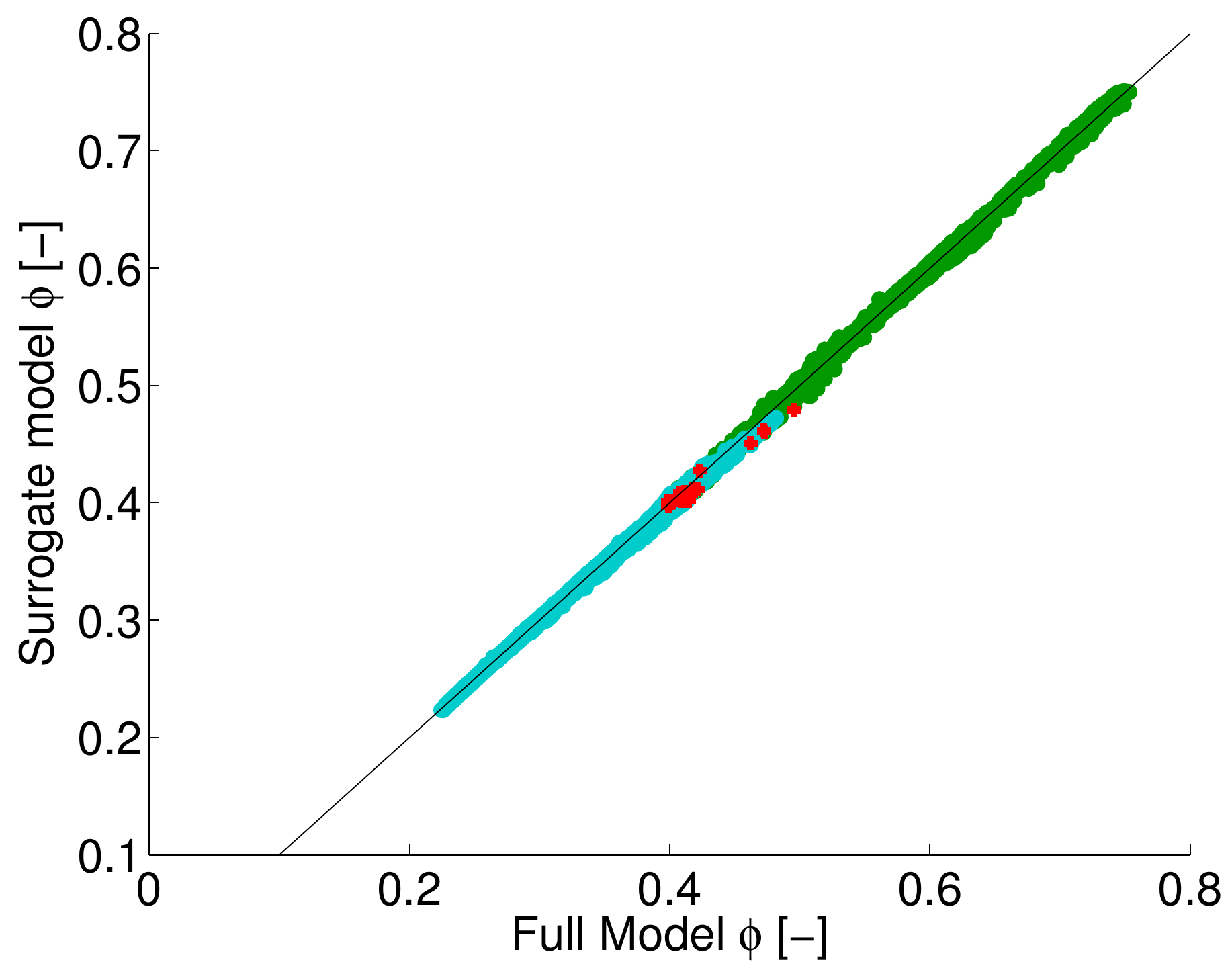} 
    \caption{scatterplot at -1350 meters}
    \label{image:scatterplot-1350}
  \end{subfigure}
  \begin{subfigure}{0.5\textwidth}
    \centering
    \includegraphics[height=5cm]{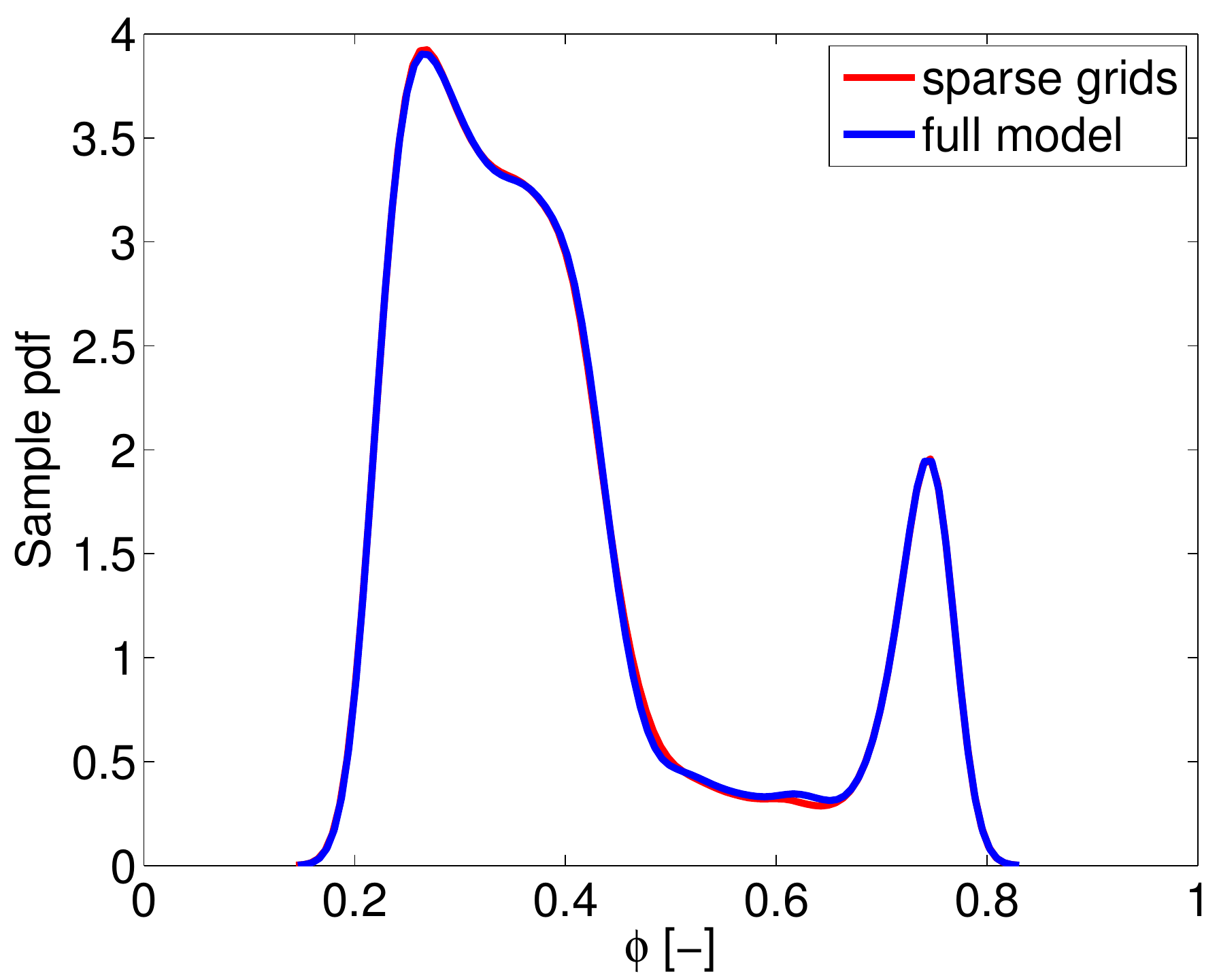} 
    \caption{pdf at -1500 meters}
    \label{image:pdf-comparison1500}
  \end{subfigure}
  \begin{subfigure}{0.5\textwidth}
    \centering
    \includegraphics[height=5cm]{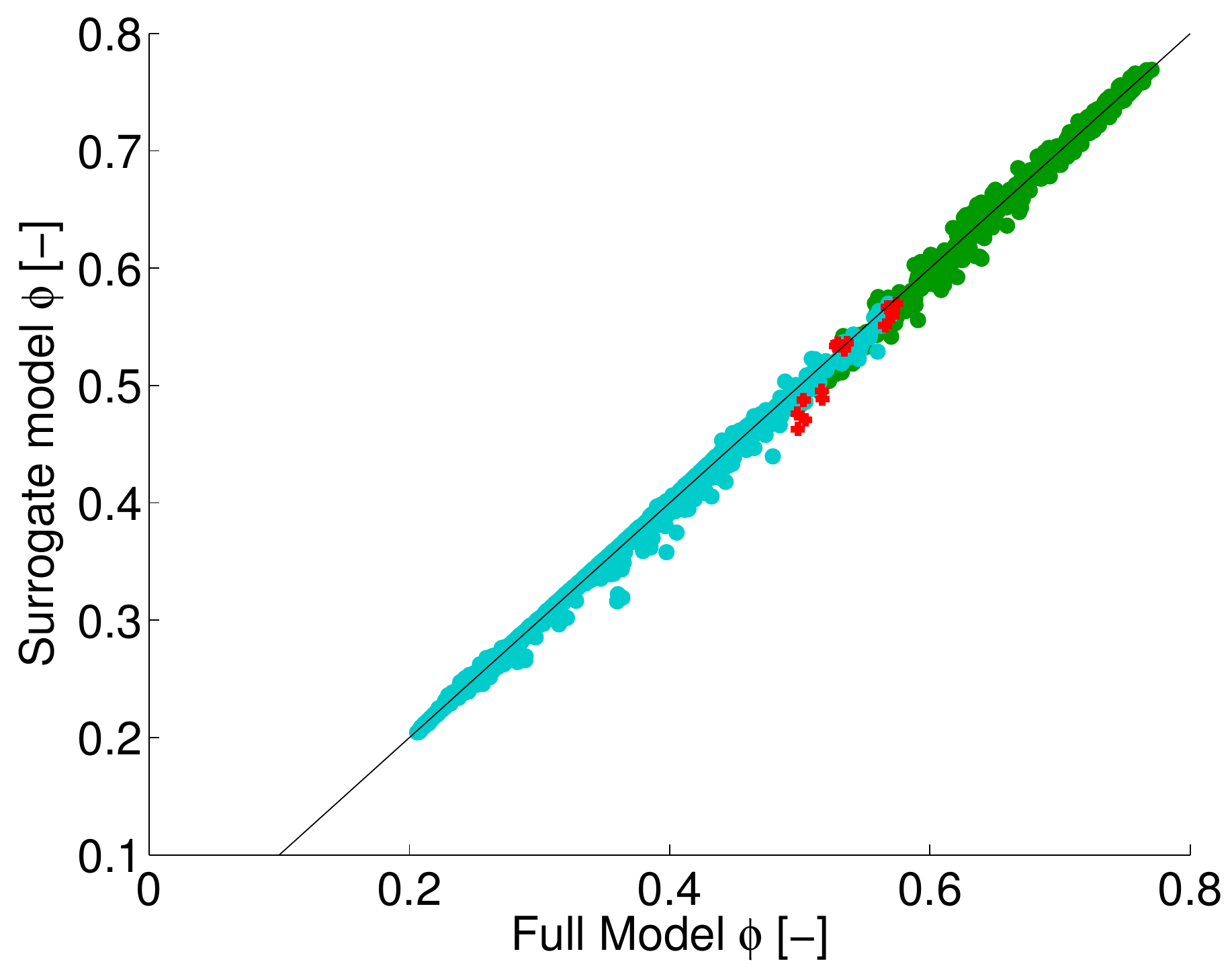}
    \caption{scatterplot at -1500 meters}
    \label{image:scatterplot-1500}
  \end{subfigure}
  \begin{subfigure}{0.5\textwidth}
    \centering
    \includegraphics[height=5cm]{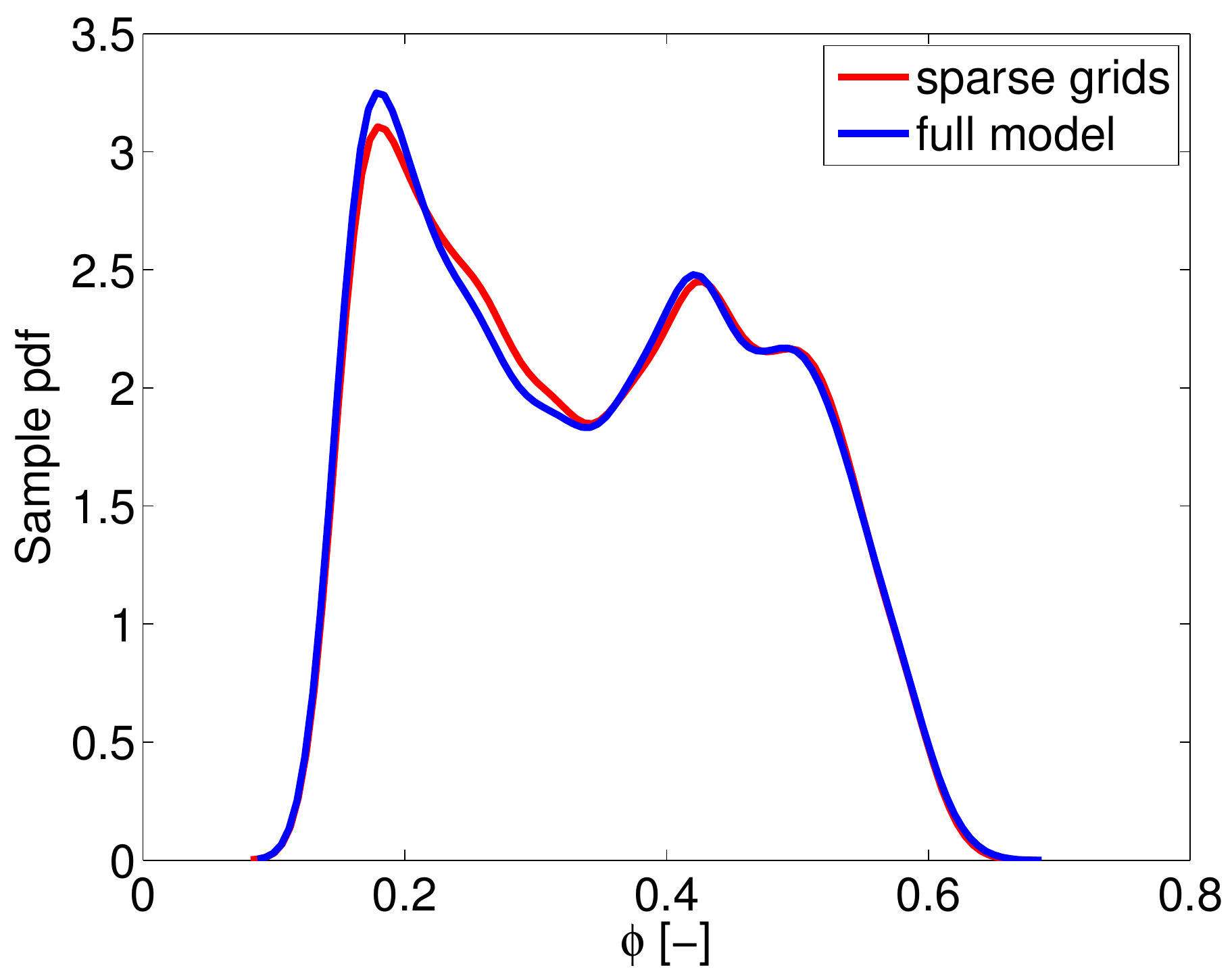} 
    \caption{pdf at -2250 meters}
    \label{image:pdf-comparison2250}
  \end{subfigure}
  \begin{subfigure}{0.5\textwidth}
    \centering
    \includegraphics[height=5cm]{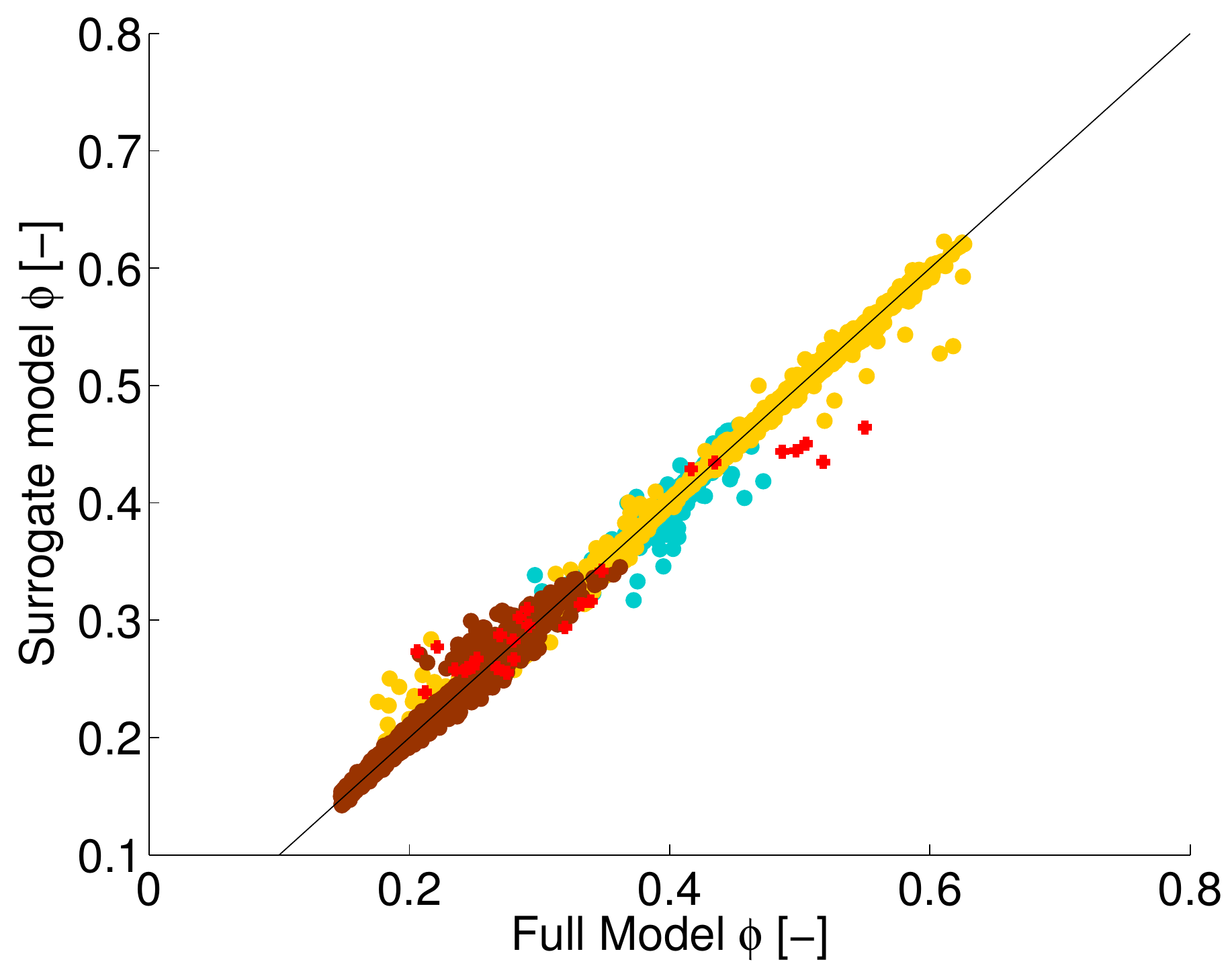}
    \caption{scatterplot at -2250 meters}
    \label{image:scatterplot-2250}
  \end{subfigure}
  \caption{Comparison of probability density functions associated to porosity obtained with full model and surrogate model at different depths.
  The color code \Lc{identifying}{identify} different materials in the scatterplots \Lc{}{and} is the same used in Figure \ref{image:material-occurrence}.}
  \label{image:pdf-comparison}
\end{figure}

\section{Conclusions}
\label{conclusions}

In this paper we devise a methodology for the quantification of uncertainty related to key outputs of sedimentary basins compaction modeling, which considers the evolution along geologic time scales of diagenesis and mechanical compaction phenomena. At these scales, information on the model parameters is typically incomplete and the availability of an uncertainty quantification procedure which is at the same time accurate from a numerical viewpoint and characterized by computational efficiency is of paramount concern.  Our methodology embeds reduced order (surrogate) modeling based on sparse grid sampling of a set of unknown parameters, which are considered to be random. 
Our work leads to the following major results:

\begin{enumerate}
	\item We extend here the approach presented in \cite{feal:compgeo,lever.eal:inversion} to a case of multi-layered domain, where we have alternating deposition of sand and shale materials. In this case, permeability can attain very low values (typically in shale materials) and exhibit variations of several order of magnitude across layers. In this context the fixed iteration method proposed in previous works (\cite{feal:compgeo}) is prone to failure. We propose here a novel one-dimensional solver grounded on a Newton iteration algorithm for the solution of the fully coupled system consisting of mass, momentum and energy conservation equations, \Lc{togheter}{together} with geochemical constitutive relations along the vertical direction. The Newton iteration algorithm clearly outperforms fixed point iterations for permeabilities within the range of those typically associated with shale-rich materials.
	\item Specific variables of the problem exhibit sharp variations moving along the vertical direction, which are located at the interface between different geomaterials. The location of such discontinuities is typically function of model parameters. The accuracy of polynomial approximations obtained from standard sparse grid methods typically  deteriorates in the presence of discontinuous mapping between input parameters and output variables. We present an original approach to tackle this issue which relies on the assumption that the depth of the interfaces among layers shows a smooth dependence on the parameters and can therefore be accurately estimated with a sparse-grids approach. We consider a synthetic test case characterized by realistic evolutionary scales. Mechanical compaction modulus and shale permeability are considered as random parameters. Sobol' indices analysis \Lc{suggests}{indicates} that vertical shale permeability has larger effect than vertical compaction modulus on the location of material interfaces where discontinuities of state variables (e.g., porosity) may occur. 
	\item We analyze the convergence of the approximation of the interface locations obtained through sparse grid sampling, based on the mean estimated location of the interfaces and on the maximum error with respect to the full model simulation within the considered random parameter space. We observe that for shallow depths  the sparse grid approximation converges rapidly when increasing the number of collocation points, independently of the strategy adopted for the sparse grid sampling approach. For deep locations the input/output mapping becomes more complex. In such conditions, the efficiency of the sparse grid approximation greatly improves when anisotropic sampling is performed, i.e. when  the order of the sparse grid approximation increases proportionally to the intensity of the Sobol' sensitivity indices associated with each parameter. The number of collocation points required to obtain a given level of accuracy may be reduced by two orders of magnitude when considering an anisotropic sparse grid instead of an isotropic one. 
	\item We finally investigate the uncertainty propagation from input parameters to model responses which exhibit a discontinuous behavior in the parameter space, i.e. porosity, considering its probability density function at selected depths. We compare the results obtained with the two-step surrogate model and the full model within a Monte Carlo framework. Results \Lc{suggest}{indicate} that the pdfs associated with the sparse grid approximation are in close agreement with those obtained with the direct simulation of the full model. We display for each selected location also a scatterplot which allows a direct comparison between the porosity values predicted with the sparse grid approximation against those given by the full model for the same location and for all considered Monte Carlo realizations of the parameters. What emerges is that even though for deeper layers the accuracy of sparse grid approximation decrease, points in the scatter plot exhibit an increasing dispersion around the bisector and the number of mismatched points increase, the entity of these inaccuracies is small enough to have only a minor consequence on the estimation of the pdf of porosity. This result suggests that the proposed methodology is suitable to be embedded within model calibration techniques, which will be considered in future works.
\end{enumerate}

\paragraph{Acknowledgement}
  F. Nobile and L. Tamellini received support from the
  Center for ADvanced MOdeling Science (CADMOS).
  L. Tamellini also received support from the 
  Gruppo Nazionale Calcolo Scientifico - Istituto Nazionale di Alta Matematica
  ``Francesco Severi'' (GNCS-INDAM).


\bibliographystyle{elsarticle-num}
\bibliography{lorenzo_biblio}

\end{document}